\documentclass[11pt]{amsart}

\usepackage{latexsym}
\usepackage{amssymb}
\usepackage{amsmath}
\usepackage{fancybox,color}
\usepackage{amssymb}
\usepackage{amsfonts}
\usepackage{latexsym}
\usepackage{amsmath,systeme}
\usepackage[leqno]{mathtools}

\usepackage{color}
\usepackage{amsxtra}
\usepackage{mathtools}
\usepackage{enumerate}
\usepackage{blindtext}
\usepackage[inline]{enumitem}
\usepackage{nicefrac}
\mathtoolsset{showonlyrefs} 
\usepackage{ mathrsfs }
\usepackage{tikz}
\tikzset{
  treenode/.style = {align=center, inner sep=0pt, text centered,
    font=\sffamily},
  arn_r/.style = {treenode, circle, blue, draw=blue,fill=yellow, 
    text width=1em, very thick},
    dot/.style={circle,draw,inner sep=1.2,fill=black},
}
\usepackage{forest}

\usetikzlibrary{shadings}
\usepackage{color}
\usepackage{hyperref}
\hypersetup{
  colorlinks   = true, 
  urlcolor     = blue, 
  linkcolor    = blue, 
  citecolor   = red 
}
\xdefinecolor{darkgreen}{RGB}{175, 193, 36}

\usepackage{amsfonts}
\usepackage{enumerate}
\usepackage{url}
\usepackage{color}      
\usepackage{epsfig}
\usepackage{graphicx}
\usepackage{esint}

\usepackage{amsmath}
\usepackage{amsthm}
\usepackage{color}
\usepackage{amssymb}
\usepackage{graphicx}
\usepackage{mathrsfs}
\usepackage{multicol}
\usepackage{wrapfig}
\usepackage{times}
\usepackage{enumerate}
\newtheorem{definition}{Definition}
\newtheorem{remark}{Remark}
\newtheorem{Lemma}{Lemma}
\newtheorem{theorem}{Theorem}

\usepackage{mathtools}
\mathtoolsset{showonlyrefs}

\def\vint_#1{\mathchoice%
          {\mathop{\kern 0.2em\vrule width 0.6em height 0.69678ex depth -0.58065ex
                  \kern -0.8em \intop}\nolimits_{\kern -0.4em#1}}%
          {\mathop{\kern 0.1em\vrule width 0.5em height 0.69678ex depth -0.60387ex
                  \kern -0.6em \intop}\nolimits_{#1}}%
          {\mathop{\kern 0.1em\vrule width 0.5em height 0.69678ex
              depth -0.60387ex
                  \kern -0.6em \intop}\nolimits_{#1}}%
          {\mathop{\kern 0.1em\vrule width 0.5em height 0.69678ex depth -0.60387ex
                  \kern -0.6em \intop}\nolimits_{#1}}}
\def\vintslides_#1{\mathchoice%
          {\mathop{\kern 0.1em\vrule width 0.5em height 0.697ex depth -0.581ex
                  \kern -0.6em \intop}\nolimits_{\kern -0.4em#1}}%
          {\mathop{\kern 0.1em\vrule width 0.3em height 0.697ex depth -0.604ex
                  \kern -0.4em \intop}\nolimits_{#1}}%
          {\mathop{\kern 0.1em\vrule width 0.3em height 0.697ex depth -0.604ex
                  \kern -0.4em \intop}\nolimits_{#1}}%
          {\mathop{\kern 0.1em\vrule width 0.3em height 0.697ex depth -0.604ex
                  \kern -0.4em \intop}\nolimits_{#1}}}

\newcommand{\aveint}[2]{\mathchoice%
          {\mathop{\kern 0.2em\vrule width 0.6em height 0.69678ex depth -0.58065ex
                  \kern -0.8em \intop}\nolimits_{\kern -0.45em#1}^{#2}}%
          {\mathop{\kern 0.1em\vrule width 0.5em height 0.69678ex depth -0.60387ex
                  \kern -0.6em \intop}\nolimits_{#1}^{#2}}%
          {\mathop{\kern 0.1em\vrule width 0.5em height 0.69678ex depth -0.60387ex
                  \kern -0.6em \intop}\nolimits_{#1}^{#2}}%
          {\mathop{\kern 0.1em\vrule width 0.5em height 0.69678ex depth -0.60387ex
                  \kern -0.6em \intop}\nolimits_{#1}^{#2}}}
\newcommand{\ol}{\overline}
\newcommand{\ul}{\underline}

\newcommand{\R}{\mathbb{R}}
\newcommand{\T}{\mathbb{T}}

\newcommand{\eps}{\varepsilon}

\newcommand{\half}{\frac{1}{2}}

\DeclareMathOperator{\med}{\,median\,}

\begin{document}
\title[Systems involving mean value formulas on trees]{Systems involving mean value formulas on trees }
\author[A. Miranda, C. A. Mosquera and J. D. Rossi]{Alfredo Miranda, Carolina A. Mosquera and Julio D. Rossi}

\address{Alfredo Miranda, Carolina A. Mosquera and Julio D. Rossi
\hfill\break\indent
Departamento  de Matem{\'a}tica, FCEyN,
Universidad de Buenos Aires,
\hfill\break\indent Pabellon I, Ciudad Universitaria (1428),
Buenos Aires, Argentina.}
\email{{\tt amiranda@dm.uba.ar, mosquera@dm.uba.ar, jrossi@dm.uba.ar}}

\begin{abstract} In this paper we study the Dirichlet problem
for systems of mean value equations on a regular tree. We deal both with
the directed case (the equations verified by the components of the system at a node in the tree
only involve values of the unknowns at the successors of the node in the tree) and the undirected case
(now the equations also involve the predecessor in the tree).
We find necessary and sufficient conditions
on the coefficients in order to have 
existence and uniqueness of solutions for continuous boundary data.
In a particular case, we also include
an interpretation of such solutions as a limit of value functions of suitable
two-players zero-sum
games. 
\end{abstract}

\date{}

\maketitle

\section{Introduction}

A tree is, informally, an infinite graph in which each node but one (the root of the tree) has exactly $m+1$ connected nodes, $m$
successors and one predecessor (see below for a precise description of a regular tree). 
Regular trees and mean value averaging operators on them play the role of being a discrete model analogous to 
the unit ball and continuous partial differential equations in it. In this sense, 
linear and nonlinear mean value properties on trees are models that are close (and related to) 
to linear and nonlinear PDEs. The literature dealing with models and equations
given by mean value formulas on trees is quire large but mainly focused on single equations. We quote  
\cite{ary,BBGS,DPMR1,DPMR2,DPMR3,KLW,KW,Ober,s-tree,s-tree1}
and references therein
for references that are closely related to our results, but
the list is far from being complete. 

Our main goal here is to look for existence and uniqueness of solutions to systems of mean value formulas 
on regular trees. When dealing with systems two main difficulties arise: the first one comes from the 
operators used to obtain the equations that govern the components of the system and the second one comes from the coupling between the components. Here we deal with linear couplings with coefficients 
in each equation that may change form one point to another and with linear or nonlinear mean value properties given in terms of averaging operators involving the successors together with a possible linear dependence on the predecessor. 

Our main result can be summarized as follows: for a general system of averaging operators with linear coupling on a regular tree we find the sharp conditions (necessary and 
sufficient conditions) on the coefficients of the coupling and the contribution of the predecessor/sucessors in such a way
that the Dirichlet problem for the system with continuous boundary data has existence and uniqueness
of solutions.

Now, let us introduce briefly some definitions and notations needed to make
precise the statements of our main results.

{\bf The ambient space, a regular tree.} 
	Given $m\in\mathbb{N}_{\ge2},$ a tree $\T$ with regular 
	$m-$branching is an infinite graph that consists of a root, denoted as
    the empty set $\emptyset$, and an infinite number of nodes, labeled as all finite sequences 
    $(a_1,a_2,\dots,a_k)$ with $k\in \mathbb{N},$ 
    whose coordinates $a_i$ are chosen from $\{0,1,\dots,m-1\}.$ 
    
    \begin{center}
        \pgfkeys{/pgf/inner sep=0.19em}
        \begin{forest}
            [$\emptyset$,
                [0
                    [0
                        [0 [,edge=dotted]]
                        [1 [,edge=dotted]]
                        [2 [,edge=dotted]]
                    ]
                    [1
                        [0 [,edge=dotted]]
                        [1 [,edge=dotted]]
                        [2 [,edge=dotted]]
                    ]
                    [2
                        [0 [,edge=dotted]]
                        [1 [,edge=dotted]]
                        [2 [,edge=dotted]]
                    ]
                ]
                [1
                    [0
                        [0 [,edge=dotted]]
                        [1 [,edge=dotted]]
                        [2 [,edge=dotted]]
                    ]
                    [1
                        [0 [,edge=dotted]]
                        [1 [,edge=dotted]]
                        [2 [,edge=dotted]]
                    ]
                    [2
                        [0 [,edge=dotted]]
                        [1 [,edge=dotted]]
                        [2 [,edge=dotted]]
                    ]
                ]
                [2
                    [0
                        [0 [,edge=dotted]]
                        [1 [,edge=dotted]]
                        [2 [,edge=dotted]]
                    ]
                    [1
                        [0 [,edge=dotted]]
                        [1 [,edge=dotted]]
                        [2 [,edge=dotted]]
                    ]
                    [2
                        [0 [,edge=dotted]]
                        [1 [,edge=dotted]]
                        [2 [,edge=dotted]]
                    ]
                ]    
            ]
        \end{forest}
        
       A tree with $3-$branching.
    \end{center}
    
    The elements in $\T$ are called vertices. 
    Each vertex $x$ has $m$ successors, obtained by adding 
    another coordinate. We will denote by 
    \[
        S(x)\coloneqq\{(x,i)\colon i\in\{0,1,\dots,m-1\}\}
    \]
    the set of successors of the vertex $x.$ 
    If $x$ is not the root then $x$ has a only an 
    immediate predecessor, which we will denote $\hat{x}.$
    The segment connecting a vertex $x$ with $\hat{x}$ is called an 
    edge and denoted by $(\hat{x},x).$
     
    A vertex $x\in\T$ has level $k\in\mathbb{N}$ if $x=(a_1,a_2,\dots,a_k)$.   
    The level of $x$ is denoted by $|x|$ and 
    the set of all $k-$level vertices is denoted by $\T^k.$
    We say that the edge $e=(\hat{x},x)$ has $k-$level if 
    $x\in \T^k.$
        
    A branch of $\T$ is an infinite sequence of vertices starting at the root, where each of the vertices in the sequence 
    is followed 
    by one of its immediate successors.
    The collection of all branches forms the boundary of $\T$, denoted 
    by $\partial\T$.
    Observe that the mapping $\psi:\partial\T\to[0,1]$ defined as
    \[
        \psi(\pi)\coloneqq\sum_{k=1}^{+\infty} \frac{a_k}{m^{k}}
    \]
    is surjective, where $\pi=(a_1,\dots, a_k,\dots)\in\partial\T$ and
    $a_k\in\{0,1,\dots,m-1\}$ for all $k\in\mathbb{N}.$ Whenever
    $x=(a_1,\dots,a_k)$ is a vertex, we set
    $
        \psi(x)\coloneqq\psi(a_1,\dots,a_k,0,\dots,0,\dots).$

{\bf Averaging operators.}
Let $F\colon\R^m\to\R$ be a continuous function. We call $F$
an averaging operator if it satisfies the following: 
$$
\begin{array}{l}
\displaystyle F(0,\dots,0)=0 \mbox{ and } F(1,\dots,1)=1; \\[5pt]
\displaystyle F(tx_1,\dots,tx_m)=t F(x_1,\dots,x_m); \\[5pt]
\displaystyle F(t+x_1,\dots,t+x_m)=t+ F(x_1,\dots,x_m),\qquad \mbox{for all $t\in\R;$} \\[5pt]
\displaystyle F(x_1,\dots,x_{m})<\max\{x_1,\dots,x_{m}\},\qquad \mbox{if not all
		   $x_j$'s are equal;} \\[5pt] 
		   \displaystyle \mbox{$F$ is nondecreasing with respect to each variable;}
		   \end{array}
		   $$
		    in addition, we will assume that
$F$ is permutation
invariant, that is,
$$
F(x_1,\dots,x_m)= F(x_{\tau(1)},\dots,x_{\tau(m)})
$$
for each permutation $\tau$ of $\{1,\dots,m\}$ and that there exists $0<\kappa<1$
such that
\begin{equation}\label{Pro.intro}
F(x_1+c,\dots,x_m)\le F(x_1,\dots,x_m)+c\kappa
\end{equation}
for all $(x_1,\dots,x_m)\in\R^m$ and for all $c>0.$

As examples of averaging operators we mention the following ones:
The first example is taken from \cite{KLW}. For $1<p<+\infty,$
	the operator
	$F^p(x_1,\dots,x_m)=t$ from $\R^m$ to $\R$  defined
	implicitly by
	\[
	\sum_{j=1}^m (x_j-t)|x_j-t|^{p-2}=0
	\]
	is a permutation invariant averaging operator.
Next, we consider, for $0\le\alpha\le1$ and $0<\beta \leq 1$ with
		$\alpha+\beta=1$
		$$
		\begin{array}{ll}
			\displaystyle F_0(x_1,\dots,x_m)=\frac{\alpha}2
			\left(\max_{1\le j\le m}
			\{x_j\}+\min_{1\le j\le m}\{x_j\}
			\right) + \frac{\beta}m\sum_{j=1}^m x_j,\\
			\displaystyle F_1(x_1,\dots,x_m)=\alpha
			\underset{{1\le j\le m}}{\med}\{x_j\}+
			 \frac{\beta}m\sum_{j=1}^m x_j,\\
			 \displaystyle F_2(x_1,\dots,x_m)=\alpha
			\underset{{1\le j\le m}}{\med}\{x_j\}+
			 \frac{\beta}2
			 \left(\max_{1\le j\le m}
			\{x_j\}+\min_{1\le j\le m}\{x_j\}
			\right),
		\end{array}
		$$
		where
		\[
\underset{{1\le j\le m}}{\med}\{x_j\}\, :=
   \begin{cases}
  y_{\frac{m+1}2} & \text{ if }m \text{ is even},  \\
  \displaystyle\frac{y_{\frac{m}2}+ y_{(\frac{m}2 +1)}}2
  & \text{ if }m \text{ is odd},
  \end{cases}
\]
with $\{y_1,\dots, y_m\}$ a nondecreasing rearrangement of
$\{x_1,\dots, x_m\}.$
$F_0, F_1$ and $F_2$ are permutation invariant
averaging operators. Note that $F_0$ and $F_1$ verify \eqref{Pro.intro} but $F_2$ does not.

Given two averaging operators $F$ and $G$ we deal with the system
\begin{equation}
\label{ED1.general}
\left\lbrace
\begin{array}{ll}
\displaystyle u(x)=(1-p_k)\Big\{ \displaystyle
(1-\beta_k^u ) F(u(y_0),...,u(y_{m-1})) + \beta_k^u u (\hat{x}) \Big\}+p_k v(x)  &  \ x\in\T_m ,  \\[10pt]
\displaystyle v(x)=(1-q_k)\Big\{ (1-\beta_k^v ) G(v(y_0),...,v(y_{m-1})) + \beta_k^v v (\hat{x}) \Big\}+q_k u(x)  &  \ x \in\T_m , 
\end{array}
\right.
\end{equation}
here $(y_i)_{i=0,...,m-1} $ are the successors of $x$. In order to have a probabilistic
interpretation of the equations in this system (see below), we will assume that $\beta_k^u $, 
$\beta_k^v$, $p_k$, $q_k$ are all in $[0,1]$ and moreover, we will also assume that $\beta_k^u$
and $\beta_k^v$ are bounded away from 1, that is, $\beta_k^u, \beta_k^v \leq c < 1$ and 
that there is no $k$ such that
$p_k=q_k=1$. 

We supplement \eqref{ED1.general} with boundary data. We take 
two continuous functions $f,g:[0,1] \mapsto \mathbb{R}$ and 
impose that along any branch of the tree we have that
\begin{equation}
\label{ED2.general}
\left\lbrace
\begin{array}{ll}
\displaystyle \lim_{{x}\rightarrow z = \psi (\pi)}u(x) = f(z) ,  \\[10pt]
\displaystyle \lim_{{x}\rightarrow z = \psi(\pi)}v(x)=g(z). 
\end{array}
\right.
\end{equation}
Here the limits are understood along the nodes in the branch as 
the level goes to infinity. That is, if the branch is given by the sequence $\pi=\{x_n\} \subset \T_m$, 
$x_{n+1} \in S(x_{n})$, then we ask for $u(x_n) \to f(\psi (\pi))$ as $n\to \infty$.

Our main result is to obtain necessary and sufficient conditions
on the coefficients $\beta_k^u$, $\beta_k^v$, $p_k$ and $q_k$
in order to have solvability of the Dirichlet problem, \eqref{ED1.general}--\eqref{ED2.general}.  

\begin{theorem} \label{teo.main.general}
For every $f,g:[0,1] \mapsto \mathbb{R}$ continuous functions
the system \eqref{ED1.general} has a unique solution satisfying \eqref{ED2.general} if and only if the 
coefficients $\beta_k^u$, $\beta_k^v$, $p_k$ and $q_k$ satisfy the following
conditions:
\begin{equation} \label{cond.general}
\begin{array}{rl}
(1) & \displaystyle \sum_{k=1}^\infty p_k < + \infty ,\\[10pt]
(2) & \displaystyle \sum_{k=1}^\infty q_k < + \infty, \\[10pt]
(3) & \displaystyle \sum_{k=1}^\infty \prod_{j=1}^k \frac{\beta_j^u}{1- \beta_j^u} < + \infty, \\[10pt]
(4) & \displaystyle \sum_{k=1}^\infty \prod_{j=1}^k \frac{\beta_j^v}{1- \beta_j^v} < + \infty .
\end{array}
\end{equation}
\end{theorem}

\begin{remark} \label{remark.cond} {\rm 
When $\beta_k^u$ is a constant $$\beta_k^u \equiv \beta$$
the condition (3) in \eqref{cond.general}
reads as
$$
\sum_{k=1}^{\infty} \prod_{j=1}^k \frac{\beta_j^u}{1-\beta_j^u}
= \sum_{k=1}^{\infty} \Big(\frac{\beta}{1-\beta} \Big)^k <+\infty
$$
and hence we get
$$
\beta < \frac12
$$
as the right condition for existence of solutions when $\beta_k^u$ is constant, $\beta_k^u \equiv \beta$.  
Also $\beta^v_j < 1/2$ is the right condition when $\beta^v_j$ is constant.

Let us briefly comment on the conditions 
\eqref{cond.general}. The third and fourth conditions imply that when $x$ is a node with $k$ large
the influence of the predecessor in the value of the components at $x$ is small (hence, there is
more influence of the successors). 
The first and second conditions say that when we are at a point $x$ with $k$ large, then the influence of the
other  component is small. Hence with these conditions one guarantees that for $x$ with large $k$
the value of the components of \eqref{ED1.general}, $u(x)$ and $v(x)$, depend 
moistly on successors of $x$ and this is exactly
what is needed to make possible to fulfill the boundary conditions \eqref{ED2.general}.
}
\end{remark}

\begin{remark} \label{NxN} {\rm
Our results can be extended to $N\times N$ systems with unknowns $(u_1,....u_N)$,
$\! u_i :\T \mapsto \mathbb{R}$ of the form
\begin{equation}
\label{ED1.general.NxN}
\left\lbrace
\begin{array}{l}
\displaystyle u_i (x)\!  = \! \Big(1 \!  - \!  \sum_{j=1}^N p_{i,j,k} \Big)\Big\{ \displaystyle
(1 \! - \!  \beta_k^i ) F_i(u(y_0),...,u(y_{m-1})) \! + \!  \beta_k^i u_i (\hat{x}) \Big\} \!  +\!  \sum_{j=1}^N p_{i,j,k} u_j (x)  , 
 \\[10pt] 
\displaystyle \lim_{{x}\rightarrow z=\psi(\pi)}u_i (x) = f_i (z). 
\end{array}
\right.
\end{equation}
Here $F_i$ is an averaging operator, $0\leq p_{i,j,k}\leq 1$ depends on the level of $x$ and on the indexes of the components 
$i,j$ and on the level of the node in the tree $k$ and are assumed to satisfy $p_{i,i,k} =0$ and 
$0\leq \sum_{j=1}^N p_{i,j,k} < 1$.
The coefficients and $0 \leq \beta_k^i <1$ depends on the level of $x$
and on the component. 
For such general systems our result says that the system \eqref{ED1.general.NxN} has a unique solution 
if and only if 
\begin{equation} \label{cond.general.NxN}
\displaystyle \sum_{k=1}^\infty \sum_{j=1}^N p_{i,j,k} < + \infty \qquad \mbox{and} \qquad 
\displaystyle \sum_{k=1}^\infty \prod_{l=1}^k \frac{\beta_l^i}{1- \beta_l^i} < + \infty, \qquad \forall i=1,...,N,
\end{equation}
hold.
}
\end{remark}

To simplify the presentation we first prove the main result, Theorem \ref{teo.main.general}, 
in the special case in which 
$\beta_k^u = \beta_k^v \equiv 0$, $p_k=q_k$, $\forall k$, and the averaging operators are given by 
\begin{equation}\label{operatores.concretos}
\begin{array}{l}
\displaystyle F (u(y_0),...,u(y_{m-1})) = \displaystyle\frac12\max_{y\in S(x)}u(y)+\frac12\min_{y\in S(x)}u(y), \\[10pt]
\displaystyle G (v(y_0),...,v(y_{m-1})) = \frac{1}{m}\sum_{y\in S(x)}v(y).
\end{array}
\end{equation}
The fact that $\beta_k^u = \beta_k^v \equiv 0$ simplifies the computations and allows us to
find an explicit solution for the special case in which the boundary data $f$ and $g$ are two constants, 
$f\equiv C_1$ and $g \equiv C_2$. 
The choice of the averaging operators in \eqref{operatores.concretos} has no special relevance but allows us to
give a game theoretical interpretation of our equations (see below). 

After dealing with this simpler case, we deal with the 
general case and prove Theorem  \ref{teo.main.general} in full generality. 
Here we have general averaging operators, $F$ and $G$, 
$\beta_k^u$ and $\beta_k^v$ can be different from zero (and are allowed
to vary depending on the level of the point $x$) and $p_k$ and $q_k$ need not be equal. 
This case is more involved and now the solution
with constant boundary data  
$f\equiv C_1$ and $g \equiv C_2$ is not explicit (but
in this case, under our conditions for existence, we will construct explicit sub and super solutions
that take the boundary values).

Our system \eqref{ED1.general} with $F$ and $G$ given by \eqref{operatores.concretos} reads as
\begin{equation}
\label{ED1.general.99}
\left\lbrace
\begin{array}{l}
\displaystyle u(x)=(1-p_k)\Big\{ \displaystyle
(1-\beta_k^u ) \Big(\frac12\max_{y\in S(x)}u(y)+\frac12\min_{y\in S(x)}u(y)\Big) + \beta_k^u u (\hat{x}) \Big\}+p_k v(x)   \\[10pt]
\displaystyle v(x)=(1-q_k)\Big\{ (1-\beta_k^v )\Big( \frac{1}{m}\sum_{y\in S(x)}v(y) \Big) + \beta_k^v v (\hat{x}) \Big\}+q_k u(x) \ 
\end{array}
\right.
\end{equation}
for $x \in\T_m$.
This system has a probabilistic interpretation that we briefly describe (see Section \ref{sect-juegos}
for more details). 
First assume that $\beta_k^u=\beta_k^v\equiv 0$.
The game is a two-player zero-sum game played in two boards (each board
is a copy of the $m-$regular tree) with the following rules: the game starts at some node 
in one of the two trees $(x_0,i)$ with $x_0\in \T$ and $i=1,2$ (we add an index to denote in which
board is the position of the game). If $x_0$ is in the first board then with probability
$p_k$ the position jumps to the other board and with probability $(1-p_k)$ the two players play
a round of a Tug-of-War game (a fair coin is tossed and the winner chooses the next
position of the game at any point among the successors of $X_0$, we refer to \cite{BRLibro}, \cite{Lewicka},
\cite{PSSW} and \cite{PS}
for more details concerning Tug-of-War games); in the second board, with probability $q_k$ the position changes to the first board, and with probability $(1-q_k)$ the position goes to one of the successors of $x_0$ with uniform 
probability. 
We take a finite level $L$ and we end the game when the position arrives to a node at that level that we call 
$x_\tau$. We also fix two 
final payoffs $f$ and $g$. This means that in the first board Player I pays to Player II the amount encoded 
by $f(\psi(x_\tau))$ while in the second board the final payoff is given by $g(\psi (x_\tau))$. Then the value function 
for this game is defined as
$$
w_L (x,i) = \inf_{S_I} \sup_{S_{II}} \mathbb{E}^{(x,i)} (\mbox{final payoff})
=  \sup_{S_{II}} \inf_{S_I} \mathbb{E}^{(x,i)} (\mbox{final payoff}).
$$
Here the $\inf$ and $\sup$ are taken among all possible strategies of the players
(the choice that the Players make at every node of what will be the next position if they play
(probability $(1-p_k)$) and 
they win the coin toss (probability $1/2$)). The final payoff is given by $f$ or $g$ according to $i_\tau =1$ or
$i_\tau=2$ (the final position of the game is in the first or in the second board).  

When $\beta_k^u$ and/or $\beta_k^v$ are not zero we add at each turn of the 
game a probability 
of passing to the predecessor of the node. 

We have that the pair of functions $(u_L,v_L)$ given by $u_L(x) = w_L (x,1)$
and $v_L (x) = w_L (x,2)$ is a solution to the system \eqref{ED1.general.99}
in the finite subgraph of the tree composed by nodes of level less than $L$.
Now we can take the limit as $L \to \infty$ in the value functions for this game and we
obtain that the limit is the unique solution to our system \eqref{ED1.general.99} that 
verifies the boundary conditions \eqref{ED2.general} in the infinite tree, see Section \ref{sect-juegos}.

\medskip

Organization of the paper: In the next section, Section \ref{sect-directed}, we deal with our system
in the special case of the directed tree, $\beta_k^u = \beta_k^v \equiv 0$ with $p_k=q_k$, and $F$ and $G$ given by \eqref{operatores.concretos}; in Section \ref{sect-undirected} we deal with the general case of two averaging
operators $F$ and $G$ and with general $\beta_k^u$, $\beta_k^v$; finally, Section \ref{sect-juegos} we include the game theoretical
interpretation of our results.

\section{A particular system on the directed tree} \label{sect-directed}

Our main goal in this section is to find necessary and sufficient conditions on the sequence of coefficients $\{p_k\}$
in order to have a solution to the system
\begin{equation}
\label{ED1}
\left\lbrace
\begin{array}{ll}
\displaystyle u(x)=(1-p_k)\Big\{ \displaystyle\frac12\max_{y\in S(x)}u(y)+\frac12\min_{y\in S(x)}u(y) \Big\}+p_k v(x) \qquad &  \ x\in\T_m ,  \\[10pt]
\displaystyle v(x)=(1-p_k)\Big(\frac{1}{m}\sum_{y\in S(x)}v(y) \Big)+p_k u(x) \qquad &  \ x \in\T_m , 
\end{array}
\right.
\end{equation}
with
\begin{equation}
\label{ED2}
\left\lbrace
\begin{array}{ll}
\displaystyle \lim_{{x}\rightarrow z=\psi(\pi)}u(x) = f(z) ,  \\[10pt]
\displaystyle \lim_{{x}\rightarrow z=\psi(\pi)}v(x)=g(z). 
\end{array}
\right.
\end{equation}
Here $f,g:[0,1]\rightarrow \R$ are continuous functions.

First, let us prove a lemma where we obtain a solution to our system when the functions $f$ and $g$ are just 
two constants.
\begin{Lemma}
	\label{Lcte} Given $C_1, C_2 \in \mathbb{R}$, 
	suppose that $$\sum_{k=0}^{\infty}p_k <\infty,$$ 
	then there exists $(u,v)$
	a solution of \eqref{ED1} and \eqref{ED2} with 
	$f\equiv C_1$ and $g\equiv C_2.$
\end{Lemma}

\begin{proof}
The solution that we are going to obtain will be constant at each level. That is, $$u(x_k)=a_k \qquad \mbox{ and }
\qquad v(x_k)=b_k$$ (here $x_k$ is any vertex at level $k$) for all $k\geq 0$. 
With this simplification the system \eqref{ED1} can be expressed as
\begin{equation}
\label{cte}
\left\lbrace
\begin{array}{ll}
\displaystyle a_k=(1-p_k)a_{k+1}+p_k b_k,  \\[5pt]
\displaystyle b_k=(1-p_k)b_{k+1}+p_k a_k, 
\end{array}
\right.
\end{equation}
for each $k\geq 0$. Then, we obtain the following system of linear equations
$$
\frac{1}{1-p_k}\begin{bmatrix}
	1 & -p_k \\
	-p_k & 1 
\end{bmatrix}\begin{bmatrix}
a_k \\
b_k
\end{bmatrix}=\begin{bmatrix}
a_{k+1} \\
b_{k+1}
\end{bmatrix}.
 $$	
 Iterating, we obtain
 \begin{equation}
 \label{Sistema}
 \left( \prod_{j=0}^{k}\frac{1}{1-p_j} \right) \prod_{j=0}^{k}\begin{bmatrix}
 1 & -p_j \\
 -p_j & 1 
 \end{bmatrix}\begin{bmatrix}
 a_0 \\
 b_0
 \end{bmatrix}=\begin{bmatrix}
 a_{k+1} \\
 b_{k+1}
 \end{bmatrix}.
 \end{equation}	
 Hence our next goal  is to analyze the convergence of the involved products as $k\to \infty$. First, we deal with
 $$
 \prod_{j=0}^{k}\frac{1}{1-p_j}.
 $$
Taking the logarithm we get
$$
\ln\Big(\prod_{j=0}^{k}\frac{1}{1-p_j}\Big)=-\sum_{j=0}^{k}\ln(1-p_j).
$$
Now, using that $\lim_{x\rightarrow 0}\frac{-\ln(1-x)}{x}=1,$ as we have $\sum_{j=0}^{\infty}p_j<\infty$ by hypothesis, we can deduce that the previous serie converges, 
$$-\sum_{j=0}^{\infty}\ln(1-p_j)=U<\infty.$$
Therefore, we have that the product also converges
$$
\prod_{j=0}^{\infty}\frac{1}{1-p_j}=e^{U}=\theta<\infty.
$$
Remark that $U>0$ and hence $1<\theta<\infty$.
 
Next let us deal with the matrices and study the convergence of
$$
\prod_{j=0}^{k}\begin{bmatrix}
 1 & -p_j \\
 -p_j & 1 
 \end{bmatrix}.
$$ 
Given $j\geq 0$, let us find the eigenvalues of 
$$
M_j=\begin{bmatrix}
	1 & -p_j \\
	-p_j & 1 
\end{bmatrix}
$$
It is easy to verify that the eigenvalues are $\{1+p_j,1-p_j\}$, and the associated eigenvectors are $[1 \ 1]$ and $[-1 \ 1]$ respectively. Is important to remark that these vectors are independent of $p_j$. Then, 
we introduce the orthogonal matrix ($Q^{-1}=Q^T$)
$$
Q=\frac{1}{\sqrt{2}}\begin{bmatrix}
1 & -1 \\
1 & 1 
\end{bmatrix},
$$  
and we have diagonalized $M_j$,
$$
M_j=Q\begin{bmatrix}
1+p_j & 0 \\
0 & 1-p_j 
\end{bmatrix}Q^T,
$$
for all $j\geq 0$. Then
$$
\prod_{j=0}^{k}M_j=Q\Big(\prod_{j=0}^{k}\begin{bmatrix}
1+p_j & 0 \\
0 & 1-p_j 
\end{bmatrix}\Big)Q^T=Q\begin{bmatrix}
\prod_{j=0}^{k}(1+p_j) & 0 \\
0 & \prod_{j=0}^{k}(1-p_j) 
\end{bmatrix}Q^T.
$$
Using similar arguments as before, 
we obtain that $$\prod_{j=0}^{\infty}(1+p_j)=\alpha<\infty \qquad \mbox{ and } \qquad 
\prod_{j=0}^{\infty}(1-p_j)=\frac{1}{\theta}=\beta.$$
Notice that $0<\beta<1$ and $1<\alpha<\infty$.
Therefore, taking the limit as $k\rightarrow\infty$ in \eqref{Sistema}, we obtain
$$
\frac{1}{\beta}Q\begin{bmatrix}
\alpha & 0 \\
0 & \beta
\end{bmatrix}Q^T\begin{bmatrix}
a_0 \\
b_0
\end{bmatrix}=\begin{bmatrix}
C_1 \\
C_2
\end{bmatrix}.
$$
Then, given two constants $C_1,C_2$ this linear system has a unique solution, $[a_0 \ b_0]$ 
(because the involved matrices are nonsingular). Once we have 
the value of $[a_0 \ b_0]$ we can obtain the values $[a_k \ b_k]$ at all levels using \eqref{Sistema}. 
The limits \eqref{ED2} are satisfied by this construction.
\end{proof}

Now, we need to introduce the following definition.
\begin{definition}
	Given $f,g: [0,1] \mapsto \mathbb{R}$ and a sequence $(p_k)_{k\geq 0}$,
	\begin{itemize}
	\item  the pair $(z,w)$ is a subsolution of \eqref{ED1} and \eqref{ED2} if
	\begin{equation}
    \left\lbrace
	\begin{array}{ll}
	\displaystyle z(x)\leq (1-p_k)\Big{\{ \displaystyle \half\max_{y\in S(x)}z(y)+\half\min_{y\in S(x)}z(y) \Big\}+p_k w(x)} \qquad &  \ x\in\T_m ,  \\[10pt]
	\displaystyle w(x)\leq (1-p_k)\Big(\frac{1}{m}\sum_{y\in S(x)}w(y) \Big)+p_k z(x) \qquad &  \ x \in\T_m , 
	\end{array}
	\right.
	\end{equation}
\begin{equation}
\left\lbrace
\begin{array}{ll}
\displaystyle \limsup_{{x}\rightarrow \psi(\pi)}z(x) \leq f(\psi(\pi)) ,  \\[10pt]
\displaystyle \limsup_{{x}\rightarrow \psi(\pi)}w(x)\leq g(\psi(\pi)). 
\end{array}
\right.
\end{equation}
\item  the pair $(z,w)$ is a supersolution of \eqref{ED1} and \eqref{ED2} if
\begin{equation}
\left\lbrace
\begin{array}{ll}
\displaystyle z(x)\geq (1-p_k)\Big\{  \displaystyle\half\max_{y\in S(x)}z(y)+\half\min_{y\in S(x)}z(y) \Big\}+p_k w(x) 
\qquad &  \ x\in\T_m ,  \\[10pt]
\displaystyle w(x)\geq (1-p_k)\Big(\frac{1}{m}\sum_{y\in S(x)}w(y) \Big)+p_k z(x) \qquad &  \ x \in\T_m , 
\end{array}
\right.
\end{equation}
\begin{equation}
\left\lbrace
\begin{array}{ll}
\displaystyle \liminf_{{x}\rightarrow \psi(\pi)}z(x) \geq f(\psi(\pi)) ,  \\[10pt]
\displaystyle \liminf_{{x}\rightarrow \psi(\pi)}w(x)\geq g(\psi(\pi)). 
\end{array}
\right.
\end{equation}
\end{itemize}	
\end{definition}

With these definitions we are ready to prove a comparison principle.

\begin{Lemma}\textbf{(Comparison Principle)}
\label{Comp}	
 Suppose that $(u_{\ast},v_{\ast})$ be a subsolution of \eqref{ED1} and \eqref{ED2}, and $(u^{\ast},v^{\ast})$ be a supersolution of \eqref{ED1} and \eqref{ED2}, then it holds that 
\begin{equation}
u_{\ast} (x)\leq u^{\ast}(x)  \quad \mbox{and} \quad v_{\ast}(x)\leq v^{\ast}(x) \qquad \forall x \in \T_m.
\end{equation}
\end{Lemma}

\begin{proof}
	Suppose, arguing by contradiction, that 
	\begin{equation}
	\max \Big\{\sup_{x\in\T_m}\{u_{\ast}-u^{\ast}\},\sup_{x\in\T_m}\{v_{\ast}-v^{\ast}\}\Big\}\geq \eta>0.
	\end{equation}
	Let 
	\begin{equation}
	Q=\Big\{ x\in\T_m : \max\{(u_{\ast}-u^{\ast})(x),(v_{\ast}-v^{\ast})(x)\}\geq\eta\Big\} \neq\emptyset.
	\end{equation}
	
	Claim \# 1: If $x\in Q$, there exists $y\in  S(x)$ such that $y\in Q$.
	\textit{Proof of the claim:} Suppose that 
	\begin{equation}
	\label{D1}
	u_{\ast}(x)-u^{\ast}(x)\geq\eta \quad \mbox{and} \quad u_{\ast}(x)-u^{\ast}(x)\geq v_{\ast}(x)-v^{\ast}(x)
	\end{equation}
	then
	\begin{equation}
	\begin{array}{l}
	u_{\ast}(x)-u^{\ast}(x)\! \leq \!
	\displaystyle (1\! - \! p_k)\Big\{ \half\max_{y\in S(x)}u_{\ast}(y)+\half\min_{y\in S(x)}u_{\ast}(y)- \half\max_{y\in S(x)}u^{\ast}(y)-\half\min_{y\in S(x)}u^{\ast}(y) \Big\}\\[10pt]
	\qquad \qquad \qquad \qquad +p_k (v_{\ast}(x)-v^{\ast}(x)).
	\end{array}
	\end{equation}
	Using \eqref{D1} in the last term we obtain
	\begin{equation}
	\begin{array}{l}
\displaystyle 
	\displaystyle  (1-p_k)u_{\ast}(x)-u^{\ast}(x) \\[10pt]
	\quad \displaystyle \leq 
	\displaystyle (1-p_k)\Big\{ \half\max_{y\in S(x)}u_{\ast}(y)+\half\min_{y\in S(x)}u_{\ast}(y)- \half\max_{y\in S(x)}u^{\ast}(y)-\half\min_{y\in S(x)}u^{\ast}(y) \Big\}.
	\end{array}
	\end{equation}
	Since $(1-p_k)$ is different from zero, using again \eqref{D1} we arrive to
		\begin{equation}
	\displaystyle  \eta\leq u_{\ast}(x)-u^{\ast}(x)\leq 
	 \displaystyle \Big(\half\max_{y\in S(x)}u_{\ast}(y)- \half\max_{y\in S(x)}u^{\ast}(y)\Big)+\Big(\half\min_{y\in S(x)}u_{\ast}(y)-\half\min_{y\in S(x)}u^{\ast}(y)\Big).
	\end{equation}
	Let $u_{\ast}(y_1)=\max_{y\in S(x)}u_{\ast}(y)$; it is clear that $u^{\ast}(y_1)\leq\max_{y\in S(x)}u^{\ast}(y)$. On the other hand, let $u^{\ast}(y_2)=\min_{y\in S(x)}u^{\ast}$; now, we have $u_{\ast}(y_2)\geq\min_{y\in S(x)}u_{\ast}(y)$. Hence
	\begin{equation}
	\eta\leq \Big(\half u_{\ast}(y_1)-\half u^{\ast}(y_1)\Big)+ \Big(\half u_{\ast}(y_2)-\half u^{\ast}(y_2)\Big).
	\end{equation}
	This implies that exists $y\in S(x)$ such that $u_{\ast}(y)-u^{\ast}(y)\geq\eta$. Thus $y\in Q$. 
	
	Now suppose the other case, that is,
	\begin{equation}
	\label{D2}
	v_{\ast}(x)-v^{\ast}(x)\geq\eta \quad \mbox{and} \quad v_{\ast}(x)-v^{\ast}(x)\geq u_{\ast}(x)-u^{\ast}(x).
	\end{equation}
	Now we use the second equation. We have
	\begin{equation}
	v_{\ast}(x)-v^{\ast}(x)\leq 
	(1-p_k) \Big(\frac{1}{m}\sum_{y\in S(x)}(v_{\ast}(y)-v^{\ast}(y)) \Big)
	+p_k (u_{\ast}(x)-u^{\ast}(x)).
	\end{equation}
	Using \eqref{D2} again we obtain
	\begin{equation}
	\eta\leq v_{\ast}(x)-v^{\ast}(x)\leq 
	\frac{1}{m}\sum_{y\in S(x)}(v_{\ast}(y)-v^{\ast}(y)).
	\end{equation}
	Hence, there exists $y\in S(x)$ such that $v_{\ast}(y)-v^{\ast}(y)\geq\eta$. Thus $y\in Q$. 
	This ends the proof of the claim.
	
	Now, given $y_0\in Q$ we have a sequence $(y_k)_{k\geq 1}$ included in a branch, $y_{k+1} \in S(y_k)$,
	with $(y_k)_{k\geq 1}\subset Q$. Hence, we have
	\begin{equation}
	u_{\ast}(y_k)-u^{\ast}(y_k)\geq\eta \quad \mbox{or} \quad v_{\ast}(y_k)-v^{\ast}(y_k)\geq\eta.
	\end{equation}
	Then, there exists a subsequence $(y_{k_j})_{k\geq 1}$ such that 
	\begin{equation} \label{oooo}
	u_{\ast}(y_{k_j})-u^{\ast}(y_{k_j})\geq\eta \quad \mbox{or} \quad v_{\ast}(y_{k_j})-v^{\ast}(y_{k_j})\geq\eta.
	\end{equation}
	
	Let us call $\lim_{k\rightarrow \infty}y_k =\pi$ the branch to which the $y_k$
	belong.
	Suppose that the first case in \eqref{oooo} holds, then
	\begin{equation}
	\liminf_{j\rightarrow\infty}u_{\ast}(x_{k_j})\leq f(\psi(\pi)) \quad \mbox{and} \quad \limsup_{j\rightarrow\infty}u^{\ast}(y_{k_j})\geq f(\psi(\pi)),
	\end{equation}
	Thus, we have
	\begin{equation}
	\begin{array}{l}
\displaystyle 
	0<\eta\leq \liminf_{j\rightarrow\infty}(u_{\ast}(y_{k_j})-u^{\ast}(y_{k_j}))\\[10pt]
	\qquad \displaystyle =\liminf_{j\rightarrow\infty}u_{\ast}(y_{k_j})-\limsup_{j\rightarrow\infty}u^{\ast}(y_{k_j})\leq f(\psi(\pi))-f(\psi(\pi))=0.
	\end{array}
	\end{equation}
	Here we arrived to a contradiction. The other case is similar. This ends the proof. 
\end{proof}

To obtain the existence of a solution to our system in
the general case ($f$ and $g$ continuous functions) we will use Perron's method.
Hence, let us introduce the following set
\begin{equation}
\label{ConjA}
\mathscr{A}=\big\{(z,w)  \colon \mbox{$(z,w)$ is a subsolution to (\ref{ED1}) -- (\ref{ED2}) } \big\}.
\end{equation}

First, we observe that $\mathscr{A}$ is not empty when $f$ and $g$ are bounded below. 

\begin{Lemma}
Given $f,g\in C([0,1])$ the set $\mathscr{A}$ verifies $\mathscr{A}\neq\emptyset$.
\end{Lemma}

\begin{proof}
	Taking $z(x_k)=w(x_k)=\min\{\min f,\min g \}$ for all $k\geq 0$, this pair is such that 
	$(z,w)\in\mathscr{A}$.
\end{proof}

Now, we prove that these functions are uniformly bounded. 

\begin{Lemma} Let $M=\max\{\max f, \max g \}$.
	If $(z,w)\in\mathscr{A}$, then $$z(x)\leq M \quad \mbox{and} \quad w(x)\leq M, 
	\qquad \forall x \in \T_m. $$
\end{Lemma}

\begin{proof}
	Suppose that the statement is false. Then there exists a vertex $x_0\in\T_m$ (in some level $k$) that $z(x_0)>M$ or $w(x_0)>M$. Suppose, in first case, that $z(x_0)>M$ and $z(x_0)\geq w(x_0)$. 
	
	Claim \# 1: There exists $y_0\in S(x_0)$ such that $z(y_0)>M$. Otherwise, 
	\begin{equation}
	\begin{array}{l}
\displaystyle 
	z(x_0)\leq (1-p_k)\Big\{\half\max_{y\in S(x_0)}z(y)+\half\min_{y\in S(x_0)}z(y)\Big\}+p_k w(x_0) \\[10pt]
	\qquad \qquad \displaystyle \leq(1-p_k)\Big\{\half\max_{y\in S(x_0)}z(y)+\half\min_{y\in S(x_0)}z(y)\Big\}+p_k z(x_0).
	\end{array}
	\end{equation}
	Then
	\begin{equation}
	(1-p_k)z(x_0)\leq (1-p_k)\Big\{\half\max_{y\in S(x_0)}z(y)+\half\min_{y\in S(x_0)}z(y)\Big\}\leq (1-p_k)\Big\{\half M+\half M\Big\}.
	\end{equation}
	From where we obtain the contradiction
	\begin{equation}
	M<z(x_0)\leq M.
	\end{equation} 
	The other case is similar.
	
	Claim \# 2: If $z(x_0)>M$ and $w(x_0)>M$ exists $y_0\in S(x_0)$ such that $z(y_0)>M$. Otherwise
		\begin{equation}
		\begin{array}{l}
\displaystyle 
	z(x_0)\leq (1-p_k)\Big\{\half\max_{y\in S(x_0)}z(y)+\half\min_{y\in S(x_0)}z(y)\Big\}+p_k w(x_0)\\[10pt]
	\qquad \qquad \displaystyle \leq(1-p_k)\Big\{\half M+\half M\Big\}+p_k M\leq M,
	\end{array}
	\end{equation}
	which is a contradiction.
	
	Claim \# 3: If $z(x_0)>M$ and $w(x_0)>M$ exists $y_0\in S(x_0)$ such that $w(y_0)>M$. Otherwise
		\begin{equation}
	w(x_0)\leq (1-p_k)\Big(\frac{1}{m}\sum_{y\in S(x_0)}w(y)\Big)+p_k x(x_0)\leq(1-p_k)M+p_k M\leq M,
	\end{equation}
which is again a contradiction.

With this three claims we can ensure that exists a (infinite) sequence $x=(x_0,x_0^1,x_0^2, \dots)$
that belongs to a branch such that $z(x_0^j)>M$ (or $w(x_0^j)>M$). Then, taking limit along this branch we obtain $$\limsup_{{x}\rightarrow \pi}z(x)>M$$ and we arrive to a contradiction. 	
\end{proof}

Now let us define
\begin{equation}
\label{uv}
u(x):=\sup_{(z,w)\in\mathscr{A}}z(x) \qquad \mbox{and} \qquad v(x):=\sup_{(z,w)\in\mathscr{A}}w(x).
\end{equation}

The next lemma proves that this pair of functions is in fact the desired solution to the system
\eqref{ED1}. 

\begin{theorem}
	\label{E1}
	The pair $(u,v)$ given by \eqref{uv} is the unique solution to \eqref{ED1}--\eqref{ED2}.
\end{theorem}

\begin{proof}
	First, let us see that $(u,v)\in\mathscr{A}$.  
	Take $(z,w)\in\mathscr{A}$ and fix $x\in\T_m$, then
	\begin{equation}
	z(x)\leq(1-p_k)\Big\{\half\max_{y\in S(x)}z(y)+\half\min_{y\in S(x)}z(y)\Big\}+p_k w(x).
	\end{equation}
	As $z\leq u$ and $w\leq v$ we obtain
	\begin{equation}
	z(x)\leq(1-p_k)\Big\{\half\max_{y\in S(x)}u(y)+\half\min_{y\in S(x)}u(y)\Big\}+p_k v(x).
	\end{equation}
	Taking supremum in the left hand side we obtain
	\begin{equation}
	u(x)\leq(1-p_k)\Big\{\half\max_{y\in S(x)}u(y)+\half\min_{y\in S(x)}u(y)\Big\}+p_k v(x).
	\end{equation}
	Analogously we obtain the corresponding inequality for $v$.
	
	On the other hand, 
	\begin{equation}
	\limsup_{{x}\rightarrow \pi}z(x)=\inf_{k\geq 0}\Big\{\sup_{j\geq k}z(x_j)\Big\}\leq f(\pi)
	\end{equation}
	for all $(z,w)\in\mathscr{A}$. Then taking supremum we obtain
	\begin{equation}
	\limsup_{{x}\rightarrow \pi}u(x)\leq f(\pi),
	\end{equation}
	and the same with $v$. Hence we conclude that $(u,v)\in\mathscr{A}$.
	
	Now, we want to see that $(u,v)$ verifies the equalities in the equations. We argue
	by contradiction. First assume that there is a point $x_0\in\T_m$ where an inequality is strict, that is,
	\begin{equation}
		u(x_0)<(1-p_k)\Big\{\half\max_{y\in S(x_0)}u(y)+\half\min_{y\in S(x_0)}u(y)\Big\}+p_k v(x_0).
	\end{equation}
Let
\begin{equation}
\delta=(1-p_k)\Big\{\half\max_{y\in S(x_0)}u(y)+\half\min_{y\in S(x_0)}u(y)\Big\}+p_k v(x)-u(x_0)>0
\end{equation}
and consider the function
\[
u_{0} (x) =\left\lbrace
\begin{array}{ll}
u(x) &  \ \ x \neq x_{0},  \\
\displaystyle u(x)+\frac{\delta}{2} &  \ \ x =x_{0} . \\
\end{array}
\right.
\]	
Observe that
\[
u_{0}(x_{0})=u(x_{0})+\frac{\delta}{2}<+(1-p_k)\Big\{\half \max_{y \in S(x_{0})}u(y) + \half \min_{y \in S(x_0)}u(y)\Big\}+p_kv(x_{0})
\]
and hence
\[
u_{0}(x_{0})<+(1-p_k)\Big\{\half \max_{y \in S(x_{0})}u_0(y) + \half \min_{y \in S(x_0)}u_0(y)\Big\}+p_kv(x_{0}).
\]
Then we have that $(u_{0},v)\in \mathscr{A}$ but $u_{0}(x_{0})>u(x_{0})$ reaching a contradiction.	

A similar argument shows that we $(u,v)$ also solves the second equation in the system.

Now we show that we attain the boundary values. 
	Given two functions $f,g\in C([0,1])$ and $\eps>0$, 
	there exists $\delta>0$ such that $|f(\psi(\pi_1))-f(\psi(\pi_2))|<\eps$ and $|g(\psi(\pi_1))-g(\psi(\pi_2))|<\eps$ if 
	$|\psi(\pi_1) - \psi(\pi_2)|<\delta$. Let  $k\in\mathbb{N}$ be such that $\frac{1}{m^k}<\delta$. We divide the interval $[0,1]$ in the $m^k$ subintervals $I_j=[\frac{j-1}{m^k},\frac{j}{m^k}]$ for $1\leq j\leq m^k$. Let us
	consider the constants
\begin{equation}
C_f^j=\min_{x\in I_j}f \quad \mbox{and} \quad C_g^j=\min_{x\in I_j}g
\end{equation}
for $1\leq j\leq m^k$. 

Now we observe that,
if we consider $\T_m^k=\{x\in\T_m : |x|=k \}$ we have $\# \T_m^k=m^k$ and given $x\in\T_m^k$ any branch that have this vertex in the $k-$th level has a limit that belongs to only one segment $I_j$. Then, we have a correspondence one to one the set $\T_m^k$ with the segments $(I_j)_{j=1}^{m^k}$. Let us call $x^j$ the vertex associated to $I_j$.

Fix now $1\leq j\leq m^k$, if we consider $x^j$ as a first vertex we obtain a tree such that the boundary (via $\phi$) is $I_j$. Using the Lemma \eqref{Lcte} in this tree we can obtain a solution of \eqref{ED1} and \eqref{ED2} with the constants $C_f^j$ and $C_g^j$.

Thus, doing this in all the vertex of $T_m^k$ we have the value of some functions so called $(\ul{u},\ul{v})$, in all vertex $x\in\T_m$ with $|x|\geq k$. Then, using the equation \eqref{ED1} we can obtain the values of the $(k-1)$-level. In fact, we have
\begin{equation}
\left\lbrace
\begin{array}{ll}
\displaystyle \ul{u}(x_{k-1})=(1-p_k)\Big\{ \half\max_{y\in S(x_{k-1})}\ul{u}(y)+\half\min_{y\in S(x_{k-1})}\ul{u}(y) \Big\}+p_k \ul{v}(x_{k-1}),  \\
\displaystyle \ul{v}(x_{k-1})=(1-p_k)\Big(\frac{1}{m}\sum_{y\in S(x_{k-1})}\ul{v}(y) \Big)+p_k \ul{u}(x_{k-1}).
\end{array}
\right.
\end{equation}
Then, if we call $$A_k=\half\max_{y\in S(x_{k-1})}\ul{u}(y)+\half\min_{y\in S(x_{k-1})}\ul{u}(y)
\quad \mbox{ and } \quad B_k=\frac{1}{m}\sum_{y\in S(x_{k-1})}\ul{v}(y)$$ we obtain the system
\begin{equation}
\frac{1}{1-p_k}\begin{bmatrix}
1 & -p_k \\
-p_k & 1 
\end{bmatrix}\begin{bmatrix}
\ul{u}(x_{k-1}) \\
\ul{v}(x_{k-1})
\end{bmatrix}=\begin{bmatrix}
A_k \\
B_k
\end{bmatrix}.
\end{equation}
Solving this system we obtain all the values at $(k-1)$-level. And so we continue until obtain values for all the tree $\T_m$. Let us observe that the pair of functions $(\ul{u},\ul{v})$ verifies
\begin{equation}
\left\lbrace
\begin{array}{ll}
\displaystyle \ul{u}(x)=(1-p_k)\Big\{ \half\max_{y\in S(x)}\ul{u}(y)+\half\min_{y\in S(x)}\ul{u}(y) \Big\}+p_k \ul{v}(x) \qquad &  \ x\in\T_m ,  \\[10pt]
\displaystyle \ul{v}(x)=(1-p_k)\Big(\frac{1}{m}\sum_{y\in S(x)}\ul{v}(y) \Big)+p_k \ul{u}(x) \qquad &  \ x \in\T_m , \\[10pt]
\end{array}
\right.
\end{equation}
\begin{equation}
\left\lbrace
\begin{array}{ll}
\displaystyle \lim_{{x}\rightarrow \pi}\ul{u}(x) = C_f^j \quad \mbox{if} \ \pi\in I_j,  \\[10pt]
\displaystyle \lim_{{x}\rightarrow \pi}\ul{v}(x)=C_g^j \quad \mbox{if} \ \pi\in I_j. 
\end{array}
\right.
\end{equation}	

We can make the same construction but using the constants
\begin{equation}
D_f^j=\max_{x\in I_j}f \quad \mbox{and} \quad D_g^j=\max_{x\in I_j}g
\end{equation}
to obtain the pair of functions $(\ol{u},\ol{v})$ that verifies 
\begin{equation}
\left\lbrace
\begin{array}{ll}
\displaystyle \ol{u}(x)=(1-p_k)\Big\{ \half\max_{y\in S(x)}\ol{u}(y)+\half\min_{y\in S(x)}\ol{u}(y) \Big\}+p_k \ol{v}(x) \qquad &  \ x\in\T_m ,  \\[10pt]
\displaystyle \ol{v}(x)=(1-p_k)\Big(\frac{1}{m}\sum_{y\in S(x)}\ol{v}(y) \Big)+p_k \ol{u}(x) \qquad &  \ x \in\T_m , \\[10pt]
\end{array}
\right.
\end{equation}
and
	\begin{equation}
	\left\lbrace
	\begin{array}{ll}
	\displaystyle \lim_{{x}\rightarrow \pi}\ol{u}(x) = D_f^j \quad \mbox{if} \ \pi\in I_j,  \\[10pt]
	\displaystyle \lim_{{x}\rightarrow \pi}\ol{v}(x)=D_g^j \quad \mbox{if} \ \pi\in I_j. \\
	\end{array}
	\right.
	\end{equation}

	Now we observe that the pair $(\ul{u},\ul{v})$ is a subsolution and $(\ol{u},\ol{v})$ is a supersolution. 
	We only need to observe that
\begin{equation}
\begin{array}{ll}
\displaystyle \lim_{{x}\rightarrow \pi}\ul{u}(x) =\limsup_{{x}\rightarrow \pi}\ul{u}(x)= C_f^j\leq f(\psi(\pi)) \quad \mbox{if} \ 
\psi(\pi)\in I_j,  \\[10pt]
\displaystyle \lim_{{x}\rightarrow \pi}\ul{v}(x)=\limsup_{{x}\rightarrow \pi}\ul{v}(x)=C_g^j\leq g(\psi(\pi)) \quad \mbox{if} \ 
\psi(\pi)\in I_j. \\
\end{array}
\end{equation}	
and
\begin{equation}
\begin{array}{ll}
\displaystyle \lim_{{x}\rightarrow \pi}\ol{u}(x) =\liminf_{{x}\rightarrow \pi}\ol{u}(x)= D_f^j\geq f(\psi(\pi)) \quad \mbox{if} \ 
\psi(\pi)\in I_j,  \\[10pt]
\displaystyle \lim_{{x}\rightarrow \pi}\ol{v}(x)=\limsup_{{x}\rightarrow \pi}\ol{v}(x)=D_g^j\geq g(\psi(\pi)) 
\quad \mbox{if} \  \psi(\pi)\in I_j. 
\end{array}
\end{equation}
By construction, $\ul{u}\leq u$ and $\ul{v}\leq v$. On the other hand, given $(z,w)\in\mathscr{A}$, by the comparison principle, we obtain $z\leq \ol{u}$ and $w\leq \ol{v}$, taking supremum, we obtain $u\leq\ol{u}$ and $v\leq\ol{v}$.

Now observe that for $\psi(\pi)\in I_j$ it holds that $C_f^j\geq f(\psi(\pi))-\eps$ and $D_f^j\leq f(\psi(\pi))+\eps$, and $C_g^j\geq g(\psi(\pi))-\eps$ and $D_g^j\leq g(\psi(\pi))+\eps$. Thus, we have
\begin{equation}
\begin{array}{l}
\displaystyle 
f(\psi(\pi))-\eps\leq C_f^j=\liminf_{{x}\rightarrow \pi}\ul{u}(x)\leq\liminf_{{x}\rightarrow \pi}u(x) \\[10pt]
\qquad \displaystyle \leq\limsup_{{x}\rightarrow \pi}u(x)\leq\limsup_{{x}\rightarrow \pi}\ol{u}(x)=D_f^j\leq f(\psi(\pi))+\eps
\end{array}
\end{equation}
and
\begin{equation}
\begin{array}{l}
\displaystyle 
g(\psi(\pi))-\eps\leq C_g^j=\liminf_{{x}\rightarrow \pi}\ul{v}(x)\leq\liminf_{{x}\rightarrow \pi}v(x)\\[10pt]
\qquad \displaystyle \leq\limsup_{{x}\rightarrow \pi}v(x)\leq\limsup_{{x}\rightarrow \pi}\ol{v}(x)=D_g^j\leq g(\psi(\pi))+\eps.
\end{array}
\end{equation}	
Hence, since $\eps$ is arbitrary we obtain
\begin{equation}
\lim_{{x}\rightarrow \pi}u(x)=f(\psi(\pi)) \quad \mbox{and} \quad \lim_{{x}\rightarrow \pi}v(x)=g(\psi(\pi))
\end{equation}
and the conclude that $(u,v)$ is a solution to \eqref{ED1} and \eqref{ED2}.

The uniqueness of solutions is a direct consequence of the comparison principle.
	Suppose that $(u_1,v_1)$ and $(u_2,v_2)$ are two solutions of \eqref{ED1} and \eqref{ED2}. Then, $(u_1,v_1)$ is a subsolution and $(u_2,v_2)$ is a supersolution. From the comparison principle we get $u_1\leq u_2$ and $v_1\leq v_2$ in $\T_m$. The reverse inequalities are obtained reversing the roles of $(u_1,v_1)$ and $(u_2,v_2)$.
\end{proof}

Next, we show that the condition $\sum_{k=0}^{\infty}p_k<\infty$ is also necessary to have existence of solution for every 
continuous boundary data. 

\begin{theorem}
	Let be $C_1\neq C_2$ two constants and $(p_k)_{k\geq 0}$ a sequence of positive numbers such that 
	$$\sum_{k=1}^{\infty}p_k=+\infty.$$ 
	Then, for any two constants $C_1$ and $C_2$ such that $C_1 \neq C_2$
	the system \eqref{ED1}--\eqref{ED2} with $f\equiv C_1$ and $g\equiv C_2$  does not have a solution.
\end{theorem}

\begin{proof}
Arguing by contradiction suppose that for every pair of constants
$C_1$, $C_2$ the system has a solution $(u,v)$. First, suppose that $C_1>C_2$. 
If we take $\ol{u}\equiv\ol{v}\equiv C_1$ this pair is a supersolution to our problem. Then, by the comparison 
lemma, we get
\begin{equation}
u(x)\leq C_1 \qquad \mbox{and} \qquad v(x)\leq C_1 \qquad \mbox{for all} \ x\in\T_m.
\end{equation}
Given a level $k\geq 0$ we have
\begin{equation}
u(x_k)=(1-p_k)\Big\{\half\max_{y\in S(x_k)}u(y)+\half\min_{y\in S(x_k)}+p_kv(x_k) \Big\}\leq (1-p_k)u(x_{k+1})+p_kv(x_k)
\end{equation}
where we have chosen $x_{k+1}$ such that
$$u(x_{k+1})=\max_{y\in S(x_k)}u(y).$$ Now using the same argument we obtain
\begin{equation}
u(x_{k+1})\leq(1-p_{k+1})u(x_{k+2})+v(x_{k+1})
\end{equation}
and hence we arrive to
\begin{eqnarray*}
u(x_k) &\leq&(1-p_k)\Big\{(1-p_{k+1})u(x_{k+2})+p_{k+1}v(x_{k+1})\Big\}+p_kv(x_k)\\
&=&(1-p_k)(1-p_{k+1})u(x_{k+2})+(1-p_k)p_{k+1}v(x_{k+1})+p_kv(x_k).
\end{eqnarray*}
Remark that the coefficients sum up to 1, that is, we have 
$$
(1-p_k)(1-p_{k+1})+(1-p_k)p_{k+1}+p_k=1.
$$
Inductively, we get
\begin{equation}
u(x_k)\leq \prod_{j=0}^{l}(1-p_{k+j})u(x_{k+j})+\sum_{j=0}^{l}a_j v(x_{k+j})
\end{equation}
where the coefficients verify
\begin{equation}
\label{suma}
\prod_{j=0}^{l}(1-p_{k+j})+\sum_{j=0}^{l}a_j=1.
\end{equation}
Now, the condition $$\sum_{j=0}^{\infty}p_j=+\infty$$ is equivalent to 
$$\prod_{j=0}^{\infty}(1-p_j)=0$$ 
(just take logarithm and use that $\ln (1+x) \sim x$ for $x\to 0$).  
Then, taking $l \to \infty$ in the equation \eqref{suma}, we get
\begin{equation}
\sum_{j=0}^{\infty} a_j=1.
\end{equation}  
But $x_k,x_{k+1},x_{k+2}, \dots $ is included in a branch in $\T$, 
then we must have 
\begin{equation} \label{lll}
\lim_{j\rightarrow\infty}v(x_{k+j})=C_2.
\end{equation}
Thus, taking the limit as $l\rightarrow\infty$ in
\begin{equation}
u(x_k)\leq \prod_{j=0}^{l}(1-p_{k+j})u(x_{k+j})+\sum_{j=0}^{l}a_j v(x_{k+j}),
\end{equation}
we obtain
\begin{equation}
u(x_k)\leq \lim_{j\rightarrow\infty}\sum_{j=0}^{\infty}a_j v(x_{k+j}).
\end{equation}
Now, if we take $k\rightarrow\infty$ and use \eqref{lll}
we obtain
\begin{equation}
C_1=\lim_{k\rightarrow\infty}u(x_k)\leq C_2
\end{equation}
which is a contradiction since we assumed that $C_1 > C_2$.

To reach a contradiction if we have $C_2 > C_1$ is analogous (we have to reverse the roles of $u$ and $v$
in the previous argument).
\end{proof}

\begin{remark} {\rm When $C_1=C_2$ if we take $u\equiv v \equiv C_1$ we have a solution to our system.
Therefore, we have proved that when $\sum_{k=1}^{\infty}p_k=+\infty$ the 
only solutions to the system \eqref{ED1}--\eqref{ED2} are the trivial ones, $u\equiv v \equiv C_1$.}
\end{remark}

\section{A general system on the undirected tree} \label{sect-undirected}


In this section we deal with the general system
\begin{equation}
\label{ED3}
\left\lbrace
\begin{array}{ll}
\displaystyle u(x)=(1-p_k)\Big\{ \displaystyle
(1-\beta_k^u ) F(u(y_0),...,u(y_{m-1})) + \beta_k^u u (\hat{x}) \Big\}+p_k v(x)  &  \ x\in\T_m ,  \\[10pt]
\displaystyle v(x)=(1-q_k)\Big\{ (1-\beta_k^v ) G(v(y_0),...,v(y_{m-1})) + \beta_k^v v (\hat{x}) \Big\}+q_k u(x)  &  \ x \in\T_m ,
\end{array}
\right.
\end{equation}
with boundary conditions
\begin{equation}
\label{ED4}
\left\lbrace
\begin{array}{ll}
\displaystyle \lim_{{x}\rightarrow \pi}u(x) = f(\psi(\pi)) ,  \\[10pt]
\displaystyle \lim_{{x}\rightarrow \pi}v(x)=g(\psi(\pi)). 
\end{array}
\right.
\end{equation}


First, we want to prove existence and uniqueness of a solution.
From the computations that we made in the previous section, we see
that the key ingredients to obtain the result are: the validity of a comparison
principle and the possibility of constructing sub and supersolutions that 
take constants as boundary values.

Now, we need to introduce the concept of sub and supersolutions for this system. 

\begin{definition}
	Given $f,g: [0,1] \mapsto \mathbb{R},$
	\begin{itemize}
	\item  the pair $(\underline{u},\underline{v})$ is subsolution of \eqref{ED3} and \eqref{ED4} if
	\begin{equation}
    \left\lbrace
\begin{array}{ll}
\displaystyle \underline{u}(x)\le(1-p_k)\Big\{(1-\beta^u_k)
F(\underline{u}(y_0),...,\underline{u}(y_{m-1}))
+\beta^u_k \underline{u}(\hat{x}) \Big\}+p_k \underline{v}(x) \qquad &  \ x\in\T_m ,  \\[10pt]
\displaystyle \underline{v}(x)\le(1-q_k)\Big\{ (1-\beta^v_k)
G(\underline{v}(y_0),...,\underline{v}(y_{m-1})) +\beta^v_k \underline{v}(\hat{x}) \Big\}+q_k \underline{u}(x)  \qquad &  \ x \in\T_m , 
\end{array}
\right.
	\end{equation}
\begin{equation}
\left\lbrace
\begin{array}{ll}
\displaystyle \limsup_{{x}\rightarrow \pi}\underline{u}(x) \leq f(\psi(\pi)) ,  \\[10pt]
\displaystyle \limsup_{{x}\rightarrow \pi}\underline{v}(x)\leq g(\psi(\pi)). 
\end{array}
\right.
\end{equation}
\item  the pair $(\overline{u},\overline{v})$ is supersolution of \eqref{ED3} and \eqref{ED4} if
	\begin{equation}
    \left\lbrace
\begin{array}{ll}
\displaystyle \overline{u}(x)\ge(1-p_k)\Big\{(1-\beta^u_k)
F(\overline{u}(y_0),...,\overline{u}(y_{m-1}))
+\beta^u_k \overline{u}(\hat{x}) \Big\}+p_k \overline{v}(x) \qquad &  \ x\in\T_m ,  \\[10pt]
\displaystyle \overline{v}(x)\ge(1-q_k)\Big\{ (1-\beta^v_k)
G(\overline{u}(y_0),...,\overline{u}(y_{m-1})) +\beta^v_k \overline{v}(\hat{x}) \Big\}+q_k \overline{u}(x)  \qquad &  \ x \in\T_m , 
\end{array}
\right.
	\end{equation}
\begin{equation}
\left\lbrace
\begin{array}{ll}
\displaystyle \liminf_{{x}\rightarrow \pi}\overline{u}(x) \ge f(\psi(\pi)) ,  \\[10pt]
\displaystyle \liminf_{{x}\rightarrow \pi}\overline{v}(x)\ge g(\psi(\pi)). 
\end{array}
\right.
\end{equation}
\end{itemize}	
\end{definition}



As before, we have a comparison principle. 

\begin{Lemma}[Comparison Principle]
	Assume that $(\ul{u},\ul{v})$ is a subsolution 
	and $(\ol{u},\ol{v})$ is a supersolution, then it holds that
	\begin{equation}
	\displaystyle \ul{u}(x)\leq \ol{u}(x)  \qquad \mbox{and} \qquad
	\displaystyle \ul{v}(x)\leq \ol{v}(x) \qquad \mbox{for all} \ x\in\T_m .
	\end{equation}
\end{Lemma}

\begin{proof} The proof starts as before,
	suppose, arguing by contradiction, that
\begin{equation}
\max\Big\{\sup_{x\in\T_m}(\ul{u}-\ol{u})(x),\sup_{x\in\T_m}(\ul{v}-\ol{v})(x)\Big\}\geq \eta>0.
\end{equation}
Let 
\begin{equation}
Q=\Big\{ x\in\T_m : \max\{(\ul{u}-\ol{u})(x),(\ul{v}-\ol{v})(x)\}\geq\eta\Big\} \neq\emptyset.
\end{equation}
Now, let us call $k_0=\min\{k : x_k\in Q\}$. 
Let $x_0\in\T_m^{k_0}$ such that $x_0\in Q$. As in the previous section
our first step is to prove the following claim:

\ul{Claim $\# 1$:} There exists a sequence $(x_0,x_1,\dots)$
inside a branch such that $x_j\in Q$ for all $j\geq 0$, $x_{j+1}\in S(x_j)$.
To prove this claim, let us begin proving that exists $y\in S(x_0)$ such that $y\in Q$. Using that $x_0\in Q$ we have to consider two cases:

\textbf{First case:} 
\begin{equation}
(\ul{u}-\ol{u})(x_0)\geq \eta \qquad \mbox{and} \qquad  (\ul{u}-\ol{u})(x_0)\geq (\ul{v}-\ol{v})(x_0).
\end{equation}
The choice of $x_0\in \T_m^{k_0}$ as a node in $Q$ that has the smallest possible level implies $(\ul{u}-\ol{u})(x_0)\geq (\ul{u}-\ol{u})(\hat{x_0})$. Then, using that $\ul{u}$ is subsolution and $\ol{u}$ is supersolution we have
\begin{equation}
\begin{array}{l}
\displaystyle (\ul{u}-\ol{u})(x_0)\leq (1-p_k)(1-\beta_k^u)\Big[
F(\underline{u}(y_0),...,\underline{u}(y_{m-1})) -
F(\ol{u}(y_0),...,\ol{u}(y_{m-1}))\Big] \\[7pt]
\qquad \qquad \qquad \displaystyle +(1-p_k) \beta_k^u(\ul{u}-\ol{u})(\hat{x_0})  +p_k (\ul{v}-\ol{v})(x_0).
\end{array}
\end{equation}
Using that $(\ul{u}-\ol{u})(x_0)\geq (\ul{v}-\ol{v})(x_0)$ we obtain
\begin{equation}
\begin{array}{l}
\displaystyle 
(\ul{u}-\ol{u})(x_0)\leq(1-\beta_k^u)
\Big[
F(\underline{u}(y_0),...,\underline{u}(y_{m-1})) -
F(\ol{u}(y_0),...,\ol{u}(y_{m-1}))\Big] \\[10pt]
\qquad \qquad \qquad \displaystyle  +\beta_k^u(\ul{u}-\ol{u})(\hat{x_0}).
\end{array}
\end{equation}
Now using that $(\ul{u}-\ol{u})(x_0)\geq (\ul{u}-\ol{u})(\hat{x_0})$ we get
\begin{equation} \label{pep}
(\ul{u}-\ol{u})(x_0)\leq \Big[F(u(\ul{u}(y_0),...,\ul{u}(y_{m-1}))- 
F(\ol{u}(y_0),...,\ol{u}(y_{m-1}))\Big].
\end{equation}
Using that $F$ is an averaging operator, this implies that there exists $y\in S(x_0)$ such that 
\begin{equation}
(\ul{u}-\ol{u})(y)\geq (\ul{u}-\ol{u})(x_0)\geq \eta.
\end{equation}
In fact, if $\max_y (\ul{u}-\ol{u})(y):=t < (\ul{u}-\ol{u})(x_0)$, using that $F$ verifies
$$F(t+x_1,\dots,t+x_m)=t+ F(x_1,\dots,x_m),$$ 
and that $F$ is nondecreasing with respect to each variable
we get 
$$
F(u(\ul{u}(y_0),...,\ul{u}(y_{m-1})) = t + F(u(\ul{u}(y_0)-t,...,\ul{u}(y_{m-1})-t)
\leq t + F(u(\ol{u}(y_0),...,\ol{u}(y_{m-1})),
$$
a contradiction with \eqref{pep}. 
We can deduce that $y\in Q$, but we also obtain $(\ul{u}-\ol{u})(y)\geq (\ul{u}-\ol{u})(x_0)$ 
a property that we are going to use later.

\textbf{Second case:} 
\begin{equation}
(\ul{v}-\ol{v})(x_0)\geq \eta \qquad \mbox{and} \qquad  (\ul{v}-\ol{v})(x_0)\geq (\ul{u}-\ol{u})(x_0).
\end{equation}
Using again that $\ul{u}$ is subsolution and $\ol{u}$ is supersolution and  we have
\begin{equation}
\begin{array}{ll}
\displaystyle (\ul{v}-\ol{v})(x_0)\leq (1-q_k)(1-\beta_k^v)
\Big[
G(\underline{v}(y_0),...,\underline{v}(y_{m-1})) -
G(\ol{v}(y_0),...,\ol{v}(y_{m-1}))\Big] \\[10pt]
\displaystyle \qquad \qquad \qquad  + (1-q_k)\beta_k^v(\ul{u}-\ol{v})(\hat{x_0}) 
+q_k (\ul{u}-\ol{u})(x_0).
\end{array}
\end{equation}
Using first $(\ul{v}-\ol{v})(x_0)\geq (\ul{u}-\ol{u})(x_0)$ and then that $(\ul{v}-\ol{v})(x_0)\geq (\ul{v}-\ol{v})(\hat{x_0})$, we  obtain
\begin{equation}
(\ul{v}-\ol{v})(x_0)\leq \Big[
G(\underline{v}(y_0),...,\underline{v}(y_{m-1})) -
G(\ol{v}(y_0),...,\ol{v}(y_{m-1}))\Big] .
\end{equation}
Arguing as before, using that $G$ is an averaging operator, 
this implies that there exists $y\in S(x_0)$ such that $y\in Q$ and
\begin{equation}
(\ul{v}-\ol{v})(y)\geq (\ul{v}-\ol{v})(x_0)\geq \eta 
\end{equation} 

Now, calling $x_1\in S(x_0)$ the node that verifies 
\begin{equation}
(\ul{u}-\ol{u})(x_1)\geq (\ul{u}-\ol{u})(x_0)\geq \eta 
\qquad \mbox{or} \qquad
(\ul{v}-\ol{v})(x_1)\geq (\ul{v}-\ol{v})(x_0)\geq \eta ,
\end{equation}
we can obtain, with the same techniques used before, a node $x_2\in S(x_1)$ such that
\begin{equation}
(\ul{u}-\ol{u})(x_2)\geq (\ul{u}-\ol{u})(x_1)\geq \eta 
\qquad \mbox{or} \qquad 
(\ul{v}-\ol{v})(x_2)\geq (\ul{v}-\ol{v})(x_1)\geq \eta. 
\end{equation}
By an inductive argument we can obtain a sequence $(x_0,x_1,x_2,\dots)\subseteq Q $ such that $x_{j+1}\in S(x_j)$. This ends the proof of the claim. 

Therefore, we can take a subsequence $(x_{j_l})_{l\geq 1}$ with the following properties
\begin{equation}
\label{x1}
(\ul{u}-\ol{u})(x_{j_l})\geq \eta \qquad \mbox{for all} \quad  l\geq 1
\end{equation} 
or
\begin{equation}
\label{x2}
(\ul{v}-\ol{v})(x_{j_l})\geq \eta \qquad \mbox{for all} \quad l\geq 1.
\end{equation} 

Suppose that \eqref{x1} is true. Let $\lim_{l\rightarrow\infty}x_{j_l}=\pi$. Then, we finally arrive to
\begin{equation}
0<\eta \leq \limsup_{l\rightarrow\infty}(\ul{u}-\ol{u})(x_{j_l})=\limsup_{l\rightarrow\infty}\ul{u}(x_{j_l})-\liminf_{l\rightarrow\infty}\ol{u}(x_{j_L})\leq f(\psi(\pi))-f(\psi(\pi))=0
\end{equation}
which is a contradiction. The argument with \eqref{x2} is simillar. This ends the proof.
\end{proof}

Now we deal with constant data on the boundary, $f\equiv C_1$ and $g\equiv C_2$.
Notice that now we only
have a supersolution to our system that takes the boundary data
(and not an explicit solution as in the previous section).

\begin{Lemma}\label{super}
	Given two constants $C_1$ and $C_2$. There exists a supersolution of \eqref{ED3} 
	such that
	$$
	\lim_{x\rightarrow \pi} u(x) = C_1
	\qquad \mbox{and} \qquad
	\lim_{x\rightarrow \pi} v(x) = C_2.
	$$
	\end{Lemma}

\begin{proof}
	We look for the desired supersolution taking 
	\begin{equation}
	u(x)=\sum_{j=k}^{\infty}r_j +C_1 \qquad \mbox{ and } \qquad 	v(x)=\sum_{j=k}^{\infty}r_j +C_2 
	\end{equation}
	for every $x\in\T_m^k$.  
	To attain the boundary conditions we need that
	$$\sum_{k=1}^{\infty}r_k<\infty.$$ 
	Indeed, is this series converge, then we have 
	$$
	\lim_{x\rightarrow \pi} u(x)= \lim_{k \to \infty} \sum_{j=k}^{\infty}r_j +C_1 = C_1
	\qquad 
	\mbox{and} \qquad
	\lim_{x\rightarrow \pi} v(x)= \lim_{k \to \infty} \sum_{j=k}^{\infty}r_j +C_2 = C_2.
	$$
	
	Now, notice that $u(x) = u(\tilde{x})$ as long as they have the same level. Therefore, using that
	$F(k,...,k) =k$, since we aim for a supersolution, from the first equation we arrive to
	\begin{equation}
	\sum_{j=k}^{\infty}r_j +C_1\geq (1-p_k)(1-\beta_k^u) \Big(\sum_{j=k+1}^{\infty}r_j+C_1\Big) +(1-p_k)\beta_k^u\Big(\sum_{j=k-1}^{\infty}r_j+C_1\Big)+p_k\Big(\sum_{j=k}^{\infty}r_j+ C_2\Big).
	\end{equation}
	We can rewrite this as
	\begin{equation}
(1-p_k)\sum_{j=k}^{\infty}r_j\geq (1-p_k)(1-\beta_k^u) \sum_{j=k+1}^{\infty}r_j +(1-p_k)\beta_k^u
\sum_{j=k-1}^{\infty}r_j + p_k(C_2-C_1).
	\end{equation}
	If we call $L=C_2-C_1$, dividing by $(1-p_k)$ we obtain
	\begin{equation}
	\sum_{j=k}^{\infty}r_j\geq (1-\beta_k^u) \sum_{j=k+1}^{\infty}r_j  +\beta_k^u
	\sum_{j=k-1}^{\infty}r_j + \frac{p_k}{(1-p_k)} L .
	\end{equation}
We can write $$\sum_{j=k}^{\infty}r_j=(1-\beta_k^u)\sum_{j=k}^{\infty}r_j+\beta_k^u\sum_{j=k}^{\infty}r_j$$ 
to obtain
\begin{equation}
(1-\beta_k^u)r_k\geq \beta_k^u r_{k-1} + \frac{p_k}{(1-p_k)} L.
\end{equation}	
Then, we have
\begin{equation}
r_k\geq \frac{\beta_k^u}{(1-\beta_k^u)} r_{k-1} + \frac{p_k}{(1-p_k)} \frac{L}{(1-\beta_k^u)}.
\end{equation}
If we iterate this inequality one more time we arrive to
\begin{equation}
r_k\geq \frac{\beta_k^u}{(1-\beta_k^u)}\Big\{ \frac{\beta_{k-1}^u}{(1-\beta_{k-1}^u)}r_{k-2}+\frac{p_{k-1}}{(1-p_{k-1})}\frac{L}{(1-\beta_{k-1}^u)} \Big\} + \frac{p_k}{(1-p_k)} \frac{L}{(1-\beta_k^u)}.
\end{equation}
Then, by an inductive argument, we obtain
\begin{equation}
\begin{array}{l}
\displaystyle 
r_k \geq \prod_{j=0}^{k-1}\frac{\beta_{k-j}^u}{(1-\beta_{k-j}^u)}r_0 + \sum_{j=1}^{k-1}
\Big(\prod_{l=0}^j \frac{\beta_{k-l}^u}{(1-\beta_{k-l}^u)}\Big)\frac{p_{k-j-1}}{(1-p_{k-j-1})}\frac{L}{(1-\beta_{k-j-1}^u)}
\\[10pt] \qquad \displaystyle +\frac{p_k}{(1-p_k)} \frac{L}{(1-\beta_k^u)}.
\end{array}
\end{equation}
Let us call $$\alpha^u_j=\prod_{l=0}^j \frac{\beta^u_{k-l}}{(1-\beta^u_{k-l})}
\qquad \mbox{ and } \qquad \ol{p}_j=\frac{p_{j+1}}{(1-p_{j+1})}\frac{L}{(1-\beta^u_{j+1})}.$$ 
We get
\begin{equation}
r_k\geq  \alpha^u_{k-1}r_0+\sum_{j=1}^{k-1}\alpha^u_j \ol{p}_{k-j}+\frac{p_k}{(1-p_k)} \frac{L}{(1-\beta_k^v)}
\end{equation}

Now, from analogous computations for the other equation in \eqref{ED3}, we obtain
\begin{equation}
\begin{array}{l}
\displaystyle 
r_k \geq \prod_{j=0}^{k-1}\frac{\beta_{k-j}^v}{(1-\beta_{k-j}^v)}r_0 + \sum_{j=1}^{k-1}\Big(\prod_{l=0}^j \frac{\beta_{k-l}^v}{(1-\beta_{k-l}^v)}\Big)\frac{q_{k-j-1}}{(1-q_{k-j-1})}\frac{L}{(1-\beta_{k-j-1}^v)}\\[10pt]
\qquad \displaystyle +\frac{q_k}{(1-q_k)} \frac{L}{(1-\beta_k^v)}
\end{array}
\end{equation}
and hence, calling $$\alpha^v_j=\prod_{l=0}^j \frac{\beta^v_{k-l}}{(1-\beta^v_{k-l})}
\qquad \mbox{ and } \qquad \ol{q}_j=\frac{q_{j+1}}{(1-q_{j+1})}\frac{L}{(1-\beta^v_{j+1})}.$$ 
We get
\begin{equation}
r_k\geq  \alpha^v_{k-1}r_0+\sum_{j=1}^{k-1}\alpha^v_j \ol{q}_{k-j}+\frac{q_k}{(1-q_k)} \frac{L}{(1-\beta_k^u)}.
\end{equation}

We can take the maximum of this two right hand sides, that is,
 \begin{equation} \label{def-rk}
 \begin{array}{rl}
r_k = & \displaystyle  \max \left\{ \alpha^u_{k-1}r_0+\sum_{j=1}^{k-1}\alpha^u_j \ol{p}_{k-j}+\frac{p_k}{(1-p_k)} \frac{L}{(1-\beta_k^v)}
\, ; \,  \right. \\[12pt]
& \displaystyle \qquad \qquad \left.
 \alpha^v_{k-1}r_0+\sum_{j=1}^{k-1}\alpha^v_j \ol{q}_{k-j}+\frac{q_k}{(1-q_k)} \frac{L}{(1-\beta_k^u)}
 \right\} .
 \end{array}
\end{equation}

Now our task is to show that 
\begin{equation}
\sum_{k=1}^{\infty}r_k<\infty.
\end{equation}

Using the hypotheses on the convergence of the series, we obtain that
$$\sum_{k=1}^{\infty}\prod_{j=0}^{k-1}\frac{\beta_{k-j}^u}{(1-\beta_{k-j}^u)}<\infty,$$  
$$\sum_{k=1}^{\infty}\alpha_k<\infty$$ and 
$$\sum_{k=1}^{\infty}\ol{p}_k<\infty,$$ since $\sum_{k=1}^{\infty}p_k<\infty$. 
Moreover $$\sum_{k=1}^{\infty}\sum_{j=1}^{k-1}\alpha_j \ol{p}_{k-j}<\infty.$$ 
This last serie is convergent since the 
Cauchy product of two sequences is convergent if both sequences are nonnegative and the respective series are convergent 
(this property is known as Merten's theorem, see \cite{Rudin}).

Similar convergences also hold for the corresponding series that involve
$\alpha^v_j$ and $\ol{q}_{k-j}$.

Then, using \eqref{def-rk} and gathering all these convergences, we conclude that
\begin{equation}
\sum_{k=1}^{\infty}r_k<\infty.
\end{equation}
This ends the proof.
\end{proof}

\begin{remark}{\rm
Notice that taking $r_0$ large we can make this supersolution as large 
as we want at the root of the tree, $\emptyset$. 
}
\end{remark}

Notice that we also have a subsolution. 

\begin{Lemma}\label{subsol}
	Given two constants $C_1$ and $C_2$. There exists a subsolution of \eqref{ED3} 
	with
	$$
	\lim_{x\rightarrow \pi} u(x) = C_1
	\qquad \mbox{and} \qquad
	\lim_{x\rightarrow \pi} v(x) = C_2.
	$$
	\end{Lemma}

\begin{proof}
Using the above lemma, we know that there exists a supersolution $(\overline{u},\overline{v})$ of \eqref{ED3} and \eqref{ED4} with $f\equiv -C_1$ and $g\equiv -C_2$.

Consider $\underline{u}=-\overline{u}$ and $\underline{v}=-\overline{v}$. Then $(\underline{u},\underline{v})$ is subsolution of \eqref{ED3} and \eqref{ED4} with $f\equiv C_1$ and $g\equiv C_2$.
\end{proof}

Now we are ready to prove existence and uniqueness of a solution when the conditions
on the coefficients, \eqref{cond.general}, hold.

\begin{theorem}\label{Thm-ppal}
Assume that the coefficients verify \eqref{cond.general}, that is,
\begin{equation}\label{condiciones}
\sum_{k=1}^{\infty} \prod_{j=1}^k \frac{\beta_j^u}{1-\beta_j^u}<+\infty, \quad 
\sum_{k=1}^{\infty}  \prod_{j=1}^k \frac{\beta_j^v}{1-\beta_j^v}<+\infty, \quad
\sum_{k=1}^{\infty}p_k<+\infty \, \mbox{ and } \, 
\sum_{k=1}^{\infty}p_k<+\infty. 
\end{equation}
Then, for every $f,g\in C([0,1])$, there exists a unique solution of \eqref{ED3} and \eqref{ED4}.
\end{theorem}

\begin{proof}
We want to prove that there exists a unique solution of \eqref{ED3} and \eqref{ED4}. 
As in the previous section, let 
\begin{equation}
\label{ConjA3}
\mathscr{A}=\Big\{(z,w)\colon (z,w)\mbox{ is subsolution of (\ref{ED3}) and (\ref{ED4})}  \Big\}.
\end{equation}

We observe that $\mathscr{A}\neq\emptyset.$ In fact, taking $z(x)=w(x)=-\max\{\|f\|_{L^{\infty}}, \|g\|_{L^{\infty}}\}=C \in\R$, we obtain that $(z,w)\in \mathscr{A}.$
Moreover, functions in $\mathscr{A}$ are bounded above. In fact, using the Comparison Principle we obtain that that $z\le C=\max\{\|f\|_{L^{\infty}}, \|g\|_{L^{\infty}}\}$ and $w\le C$ for all $(z,w)\in \mathscr{A}.$

We let
\[
(u(x), v(x))= \sup_{(z,w)\in\mathscr{A}} (z,w).
\]
We want to prove that this pair of functions is solution of \eqref{ED3} and \eqref{ED4}.

If $(z,w)\in\mathscr{A},$ 
\begin{align*}
z(x)&\le(1-p_k)\Big\{ (1-\beta_k^u)\Big(\displaystyle
F(z(y_0),...,z(y_{m-1})) +\beta_k^u z(\hat{x}) \Big\}+p_k w(x)\\
&\le(1-p_k)\Big\{ (1-\beta_k^u)\Big(\displaystyle
F(u(y_0),...,u(y_{m-1}))\Big) +\beta_k^u z(\hat{x}) \Big\}+p_k v(x)
\end{align*}
and
\begin{align*}
w(x)&\le(1-q_k)\Big\{(1-\beta_k^v)G(w(y_0),...,w(y_{m-1}))+\beta_k^v w(\hat{x}) \Big\}+q_k z(x)\\
&\le(1-q_k)\Big\{ (1-\beta_k^v)G(v(y_0),...,v(y_{m-1}))+\beta_k^v v(\hat{x}) \Big\}+q_k u(x).
\end{align*}
Then, taking supremum in the left hand sides, we obtain that 
$(u,v)$ is a subsolution of the equations in \eqref{ED3}. 

For the two functions $f,g\in {C}([0,1])$ and $\eps>0$, there exists $\delta>0$ such that $|f(\psi(\pi_1))-f(\psi(\pi_2))|<\eps$ and $|g(\psi(\pi_1))-g(\psi(\pi_2))|<\eps$ if $|\psi(\pi_1) - \psi(\pi_2)|<\delta$. Let us take $k\in\mathbb{N}$ such that $\frac{1}{m^k}<\delta$. We divide the segment $[0,1]$ in $m^k$ subsegments $I_j=[\frac{j-1}{m^k},\frac{j}{m^k}]$ for $1\leq j\leq m^k$. Let us consider the constants
\begin{equation}
C_1^j=\max_{x\in I_j}f \quad , \quad C_2^j=\max_{x\in I_j}g
\end{equation}
for $1\leq j\leq m^k$. 

If we consider $\T_m^k=\{x\in\T_m : |x|=k \}$ we have $\# \T_m^k=m^k$ and given $x\in\T_m^k$ any branch that have this vertex in the k-level ends in only one segment $I_j$. Then we can relate one to one the set $\T_m^k$ with the segments $(I_j)_{j=1}^{m^k}$. Let us call $x^j$ the vertex associated to $I_j$.

Fix now $1\leq j\leq m^k$, if we consider $x^j$ as a first vertex we obtain a tree such that the boundary (via $\phi$) is $I_j$. Using the Lemma \ref{super} in this tree we can obtain a supersolution $(\overline{u},\overline{v})$ of \eqref{ED3} and \eqref{ED4} with the constants $C_1^j$ and $C_2^j$ such that $\overline{u}(x_k)=C+1$ and $\overline{v}(x_k)=C+1,$ where $C=\max\{\|f\|_{L^{\infty}}, \|g\|_{L^{\infty}}\}.$

Let 
\begin{equation}
\overline{\overline{u}}(x)=\left\lbrace
\begin{array}{ll}
\displaystyle \overline{u}(x), \quad &x\in\T_m^j \\[7pt]
\displaystyle C+1, \quad &\mbox{ otherwise } 
\end{array}
\right.
\end{equation}
and
\begin{equation}
\overline{\overline{v}}(x)=\left\lbrace
\begin{array}{ll}
\displaystyle \overline{v}(x), \quad &x\in\T_m^j \\[7pt]
\displaystyle C+1, \quad &\mbox{ otherwise. } 
\end{array}
\right.
\end{equation}
Then $(\overline{\overline{u}},\overline{\overline{v}})$ is a supersolution of \eqref{ED3} and \eqref{ED4} in $\T_m.$
Using the Comparison Principle, we get that 
\[
z(x)\le \overline{\overline{u}}(x)\, \mbox{ and } \,w(x)\le \overline{\overline{v}}(x), \qquad \mbox{for every } (z,w)\in\mathscr{A}.
\]
Then, taking supremums in the right hand sides, we conclude that
\[
u(x)\le \overline{\overline{u}}(x)\, \mbox{ and } \,v(x)\le \overline{\overline{v}}(x).
\]
Hence, we conclude that 
\[
\limsup_{x\rightarrow \pi} u(x) \leq f(\psi(\pi)).
\]
Similarly, 
\[
\limsup_{x\rightarrow \pi} v(x) \leq g(\psi(\pi)).
\]
Then $(u,v)\in\mathscr{A}.$

We want to prove that $(u,v)$ satisfies \eqref{ED3}. We know that it is a subsolution. Suppose that there exits $x_0$ such that 
\begin{equation} \label{ppp}
u(x_0)<(1-p_k)\Big\{ (1-\beta_k^u)
F(u(y_0),...,u(y_{m-1}))
 +\beta_k^u u(\hat{x_0}) \Big\}+p_k v(x_0).
\end{equation}
Let 
\begin{equation}
u^*(x)=\left\lbrace
\begin{array}{ll}
\displaystyle u(x)+\eta, \quad &x=x_0 \\[5pt]
\displaystyle u(x), \quad &x\neq x_0. 
\end{array}
\right.
\end{equation}
Since we have a strict inequality in \eqref{ppp} and $F$ is monotone and continuous, 
it is easy to check that for $\eta$ small 
$(u^*,v)\in\mathscr{A}.$ This is a contradiction because we have 
$$u^* (x_0)>u (x_0)=\sup_{(w,z)\in\mathscr{A}}w(x_0).$$ 
A similar argument shows that $(u,v)$ also solves the second equation in \eqref{ED3}.

Up to this point we have that $(u,v)$ satisfies \eqref{ED3} together with $$\limsup_{x\rightarrow \pi} u(x) \leq f(\pi) \qquad \mbox{ and } \qquad 
\limsup_{x\rightarrow \pi} v(x) \leq g(\pi).$$
Hence, our next fast is to prove that $(u,v)$ satisfies the limits in \eqref{ED4}.

By Lemma \ref{subsol}, there exists $(\tilde{u},\tilde{v})$ a subsolution of \eqref{ED3} and \eqref{ED4} in $\T_m^j,$ with
$\tilde{u}(x_j)=-(C+1)$ and $\tilde{v}(x_j)=-(C+1).$
Let 
\begin{equation}
\tilde{\tilde{u}}(x)=\left\lbrace
\begin{array}{ll}
\displaystyle \tilde{u}(x), \quad &x\in\T_m^j \\[7pt]
\displaystyle -(C+1), \quad &\mbox{ otherwise } 
\end{array}
\right.
\end{equation}
and
\begin{equation}
\tilde{\tilde{v}}(x)=\left\lbrace
\begin{array}{ll}
\displaystyle \tilde{v}(x), \quad &x\in\T_m^j \\[7pt]
\displaystyle -(C+1), \quad &\mbox{ otherwise. } 
\end{array}
\right.
\end{equation}
Then $(\tilde{\tilde{u}},\tilde{\tilde{v}})$ is a subsolution of \eqref{ED3} and \eqref{ED4} in $\T_m.$
From the definition of $(u,v)$ a the supremum of subsolutions we get
\[
u(x)\ge \tilde{\tilde{u}}(x)\, \mbox{ and } \,v(x)\ge \tilde{\tilde{v}}(x).
\]
Then we get that  $$\liminf_{x\rightarrow \pi} u(x) \ge f(\psi(\pi)) \qquad \mbox{and} \qquad 
\liminf_{x\rightarrow \pi} v(x) \ge g(\psi(\pi))$$
and hence we conclude \eqref{ED4},
$$\lim_{x\rightarrow \pi} u(x)= f(\psi(\pi)) \qquad \mbox{and} 
\qquad \lim_{x\rightarrow \pi} v(x)= g(\psi(\pi)).$$

Finally, we observe that the comparison principle gives us uniqueness of solutions. 
\end{proof} 

Finally, the nonexistence of solutions when one of the conditions fails
completes the result.

\begin{theorem}
Let be $C_1\neq C_2$ two constants and suppose that one of the following conditions
\small{\begin{align}
&\sum_{n=1}^{+\infty} \left[\prod_{j=1}^n\left(\frac{\beta_j^u}{1-\beta_j^u}\right)\right]<+\infty, 
\qquad \sum_{n=1}^{+\infty} \left[\prod_{j=1}^n\left(\frac{\beta_j^v}{1-\beta_j^v}\right)\right]<+\infty\\
&\sum_{n=1}^{+\infty}\sum_{j=1}^{n}\left(\prod_{l=j+1}^n \frac{\beta_l^u}{1-\beta_l^u}\right)\left(\frac{p_j}{1-p_j}\right)\left(\frac{1}{1-\beta_j^u}\right)<+\infty,  \\
& \sum_{n=1}^{+\infty}\sum_{j=1}^{n}\left(\prod_{l=j+1}^n \frac{\beta_l^v}{1-\beta_l^v}\right)\left(\frac{q_j}{1-q_j}\right)\left(\frac{1}{1-\beta_j^v}\right)<+\infty
\end{align}}
is not satisfied. Then, there exist two constants $C_1$ and $C_2$
such that the system \eqref{ED3} with condition \eqref{ED4} and $f\equiv C_1$ and $g\equiv C_2$ does
not have a solution.
\end{theorem}

\begin{proof}
Suppose that the system have a solution $(u,v)$ with boundary condition \eqref{ED4} with $f\equiv C_1$ and $g\equiv C_2.$ Take $C_1>C_2.$
We have
\begin{equation}
\label{ED5}
\left\lbrace
\begin{array}{ll}
\displaystyle u(x)=(1-p_k)\Big\{ (1-\beta_k^u)
F(u(y_0),...,u(y_{m-1})) +\beta_k^u u(\hat{x}) \Big\}+p_k v(x)  \qquad &  \ x \in\T_m^k , \\[10pt]
\displaystyle v(x)=(1-q_k)\Big\{(1-\beta_k^v)
G(v(y_0),...,v(y_{m-1}))+\beta_k^v v(\hat{x}) \Big\} +q_k u(x) \qquad &  \ x\in\T_m^k .
\end{array}
\right.
\end{equation}

Let us follow the path given by the maximums among successors, that is, we let 
$$u(x_1)=\max_{y\in S(\emptyset)}u(y) \qquad \mbox{and} \qquad u(x_n)=\max_{y\in S(x_{n-1})}u(y).$$ 
Then, using that $F(z_1,...,z_m) \leq \max_i z_i$, we have
\[
u(x_n)\le (1-p_n)(1-\beta_n^u) u(x_{n+1})+ (1-p_n)\beta_n^u u(x_{n-1})+p_n v(x_n),
\]
that is, 
\[
\begin{array}{l}
\displaystyle 
0 \le (1-p_n)(1-\beta_n^u) (u(x_{n+1})-u(x_n)) \\[10pt]
\displaystyle 
\qquad \qquad + (1-p_n)\beta_n^u (u(x_{n-1})-u(x_n))+p_n (v(x_n)-u(x_n)).
\end{array}
\]
If we call $a_n=u(x_{n+1})-u(x_n)$ we get
\[
0 \le (1-p_n)(1-\beta_n^u) a_n - (1-p_n)\beta_n^u a_{n-1}+p_n (v(x_n)-u(x_n)).
\]
Then
\[
0 \le (1-\beta_n^u) a_n - \beta_n^u a_{n-1}+\frac{p_n}{(1-p_n)} (v(x_n)-u(x_n)).
\]
Now, calling $b_n=v(x_n)-u(x_n)$, we obtain
\[
a_n\ge \frac{\beta_n^u}{(1-\beta_n^u)} a_{n-1} +\frac{p_n}{(1-p_n)(1-\beta_n^u)}(-b_n).
\]
Now, using the same argument one more time, we get
\[
\begin{array}{l}
\displaystyle 
a_n\ge \frac{\beta_n^u}{(1-\beta_n^u)}\frac{\beta_{n-1}^u}{(1-\beta_{n-1}^u)} a_{n-2} \\[10pt]
\displaystyle \qquad \qquad +\frac{\beta_n^u}{(1-\beta_n^u)}\frac{p_{n-1}}{(1-p_{n-1})(1-\beta_{n-1}^u)}(-b_{n-1})+ \frac{p_{n}}{(1-p_{n})(1-\beta_{n}^u)}(-b_{n}).
\end{array}
\]
Inductively, we arrive to
\[
a_n\ge \left[\prod_{j=k}^n\left(\frac{\beta_j^u}{1-\beta_j^u}\right)\right]a_k+ \sum_{j=k}^{n}\left(\prod_{l=j+1}^n \frac{\beta_l^u}{1-\beta_l^u}\right)\left(\frac{p_j}{1-p_j}\right)\left(\frac{1}{1-\beta_j}\right)(-b_j),
\]
that taking $k=1$ is
\[
a_n\ge \left[\prod_{j=1}^n\left(\frac{\beta_j^u}{1-\beta_j^u}\right)\right]a_1+ \sum_{j=1}^{n}\left(\prod_{l=j+1}^n \frac{\beta_l^u}{1-\beta_l^u}\right)\left(\frac{p_j}{1-p_j}\right)\left(\frac{1}{1-\beta_j}\right)(-b_j).
\]
On the other hand, we have that 
\[
\sum_{n=1}^M a_n= u(x_{M+1})-u(\emptyset).
\]
Then, we obtain
\begin{equation}\label{nosol}
\begin{array}{l}
\displaystyle 
u(x_{M+1})\geq u(\emptyset) +\sum_{n=1}^M \left[\prod_{j=1}^n\left(\frac{\beta_j^u}{1-\beta_j^u}\right)\right]a_1\\[10pt]
\displaystyle \qquad \qquad \qquad
+ \sum_{n=1}^M\sum_{j=1}^{n}\left(\prod_{l=j+1}^n \frac{\beta_l^u}{1-\beta_l^u}\right)\left(\frac{p_j}{1-p_j}\right)\left(\frac{1}{1-\beta_j}\right)(-b_j).
\end{array}
\end{equation}

We observe that the boundary conditions $$\lim_{j\to+\infty}u(x_j)=C_1\qquad \mbox{ and }\qquad 
\lim_{j\to+\infty}v(x_j)=C_2,$$ implies that $$\lim_{j\to+\infty}b_j=C_2-C_1.$$ 
Therefore, since we have taken $C_1>C_2$, there exists a constant $c$ such that $$(-b_j)\geq c>0$$ for $j$ large enough.

If 
\[
\sum_{n=1}^{+\infty} \left[\prod_{j=1}^n\left(\frac{\beta_j^u}{1-\beta_j^u}\right)\right]=+\infty,
\]
we obtain a contradiction from \eqref{nosol} taking the limit as $M\to \infty$. Notice that we can always 
assume that $a_1\neq 0$ (otherwise, start the argument at the first level such that $a_i \neq 0$).

Now, if 
\begin{equation}\label{rg}
\sum_{n=1}^{+\infty} \left[\prod_{j=1}^n\left(\frac{\beta_j^u}{1-\beta_j^u}\right)\right]<+\infty,
\end{equation}
but
\begin{equation} \label{www}
\sum_{n=1}^{+\infty}\sum_{j=1}^{n}\left(\prod_{l=j+1}^n \frac{\beta_l^u}{1-\beta_l^u}\right)\left(\frac{p_j}{1-p_j}\right)\left(\frac{1}{1-\beta_j^u}\right)=+\infty,
\end{equation}
we obtain again a contradiction from \eqref{nosol} 
using that $(-b_j)\geq c>0$ for $j$ large enough and taking the limit as $M\to \infty$.
Remark that when \eqref{rg} holds \eqref{www} is equivalent to 
$$
\sum_{j=1}^{+\infty} p_j = +\infty.
$$

Similarly, we can arrive to a contradiction from 
\[
\sum_{n=1}^{+\infty} \left[\prod_{j=1}^n\left(\frac{\beta_j^v}{1-\beta_j^v}\right)\right]=+\infty,
\]
or
\[
\sum_{n=1}^{+\infty}\sum_{j=1}^{n}\left(\prod_{l=j+1}^n \frac{\beta_l^v}{1-\beta_l^v}\right)\left(\frac{q_j}{1-q_j}\right)\left(\frac{1}{1-\beta_j^v}\right)=+\infty,
\]
using the second equation in \eqref{ED3} (in this case we follow a path formed by maximums 
among values on successors of the second component of the system, $v$, and start with $C_1<C_2$).
\end{proof}

\section{Game theoretical interpretation} \label{sect-juegos}

Recall that in the introduction we mentioned that the 
system \eqref{ED1.general} with $F$ and $G$ given by \eqref{operatores.concretos},
\begin{equation}
\label{ED1.general.9977}
\left\lbrace
\begin{array}{l}
\displaystyle u(x)=(1-p_k)\Big\{ \displaystyle
(1-\beta_k^u ) \Big(\frac12\max_{y\in S(x)}u(y)+\frac12\min_{y\in S(x)}u(y)\Big) + \beta_k^u u (\hat{x}) \Big\}+p_k v(x),  \\[10pt]
\displaystyle v(x)=(1-q_k)\Big\{ (1-\beta_k^v )\Big( \frac{1}{m}\sum_{y\in S(x)}v(y) \Big) + \beta_k^v v (\hat{x}) \Big\}+q_k u(x), 
\end{array}
\right.
\end{equation}
for $x \in\T_m$,
has a probabilistic interpretation. In this final section we present the details. 

The game is a two-player zero-sum game played in two boards (each board
is a copy of the $m-$regular tree) with the following rules: the game starts at some node 
in one of the two trees $(x_0,i)$ with $x_0\in \T$ and $i=1,2$ (we add an index to denote in which
board is the position of the game). If $x_0$ is in the first board then with probability
$p_k$ the position jumps to the other board and with probability $(1-p_k)(1-\beta_k^u)$ the two players play
a round of a Tug-of-War game (a fair coin is tossed and the winner chooses the next
position of the game at any point among the successors of $X_0$, we refer to \cite{BRLibro}, \cite{Lewicka},
\cite{PSSW} and \cite{PS}
for more details concerning Tug-of-War games) and with probability $(1-p_k)\beta_k^u$
the position of the game goes to the predecessor (in the first board); 
in the second board, with probability $q_k$ the position changes to the first board, and with probability $(1-q_k)
(1-\beta_k^v)$ the position goes to one of the successors of $x_0$ with uniform 
probability while with probability $(1-q_k)\beta_k^v$ then position goes to the predecessor. 
We take a finite level $L$ (large) and we add the rule that the game ends 
when the position arrives to a node at level $L$, $x_\tau$. We also have two 
final payoffs $f$ and $g$. This means that in the first board Player I pays to Player II the amount encoded 
by $f(\psi(x_\tau))$ while in the second board the final payoff is given by $g(\psi (x_\tau))$. Then the value function 
for this game is defined as
$$
w_L (x,i) = \inf_{S_I} \sup_{S_{II}} \mathbb{E}^{(x,i)} (\mbox{final payoff})
=  \sup_{S_{II}} \inf_{S_I} \mathbb{E}^{(x,i)} (\mbox{final payoff}).
$$
Here the $\inf$ and $\sup$ are taken among all possible strategies of the players
(the choice that the Players make at every node of what will be the next position if they play
(probability $(1-p_k)(1-\beta_k^u)$) and 
they win the coin toss (probability $1/2$)). The final payoff is given by $f$ or $g$ according to $i_\tau =1$ or
$i_\tau=2$ (the final position of the game is in the first or in the second board).  
The value of the game $w_L (x,i)$ encodes the amount that the players expect to get/pay
playing their best with final payoffs $f$ and $g$ at level $L$.

We have that the pair of functions $(u_L,v_L)$ given by $u_L(x) = w_L (x,1)$
and $v_L (x) = w_L (x,2)$ is a solution to the system \eqref{ED1.general.9977}
in the finite subgraph of the tree composed by nodes of level less than $L$.

Notice that the first equation encodes all the possibilities for the next position of the game in the first board. We have
\begin{eqnarray*}
\displaystyle u(x)&=&(1-p_k)\Big\{ \displaystyle
(1-\beta_k^u ) \Big(\frac12\max_{y\in S(x)}u(y)+\frac12\min_{y\in S(x)}u(y)\Big) + \beta_k^u u (\hat{x}) \Big\}+p_k v(x)
\\
 \displaystyle &=&(1-p_k)
(1-\beta_k^u ) \Big(\frac12\max_{y\in S(x)}u(y)+\frac12\min_{y\in S(x)}u(y)\Big) + (1-p_k) \beta_k^u u (\hat{x}) 
+p_k v(x).
\end{eqnarray*}
Now, we observe that the value of the game at one node $x$ in the first board is the sum of the conditional expectations:
the probability of playing $(1-p_k)
(1-\beta_k^u )$ times the value of one round of Tug-of-War (with probability $1/2$ the first player chooses
the successor that maximizes $u(y)$ and with probability $1/2$ the other player chooses $y$ such that the minimum of $u$ is achieved); plus, the probability $(1-p_k)\beta_k^u$ times the value of $u$ at the predecessor;
plus, finally, the probability of jumping to the other board, $p_k$ times the value of the 
game if this happens, $v(x)$. 

Similarly, the second equation
$$
v(x)=(1-q_k)\Big\{ (1-\beta_k^v )\Big( \frac{1}{m}\sum_{y\in S(x)}v(y) \Big) + \beta_k^v v (\hat{x}) \Big\}+q_k u(x)
$$
takes into account all the possibilities for the game in the second board. 

Remark that when $\beta_k^u=1$ then at a node of level $k$ 
there is no possibility to go to a successor (when the players play the
only possibility is to go to the predecessor). Therefore, when  $\beta_k^u=\beta_k^v=1$
this game is not well defined (since for $L$ larger than $k$ the game never ends).
Therefore our assumption that $\beta_k^u$ $\beta_k^v$ are uniformly bounded away from 1 seems
reasonable. Notice that the game is also not well defined when $p_k=q_k=1$.

Now our goal is to take the limit as $L \to \infty$ in these value functions for this game and
obtain that the limit is the unique solution to our system \eqref{ED1.general.9977} that 
verifies the boundary conditions 
\begin{equation}
\label{ED2.general.9977}
\left\lbrace
\begin{array}{ll}
\displaystyle \lim_{{x}\rightarrow z = \psi (\pi)}u(x) = f(z) ,  \\[10pt]
\displaystyle \lim_{{x}\rightarrow z = \psi(\pi)}v(x)=g(z). 
\end{array}
\right.
\end{equation}

\begin{theorem} Fix two continuous functions $f,g:[0,1] \to \mathbb{R}$. 
The values of the game $(u_L,v_L)$, that is, the solutions to \eqref{ED1.general.9977}
in the finite subgraph of the tree with nodes of level less than $L$ and conditions $u_L(x) = f (\psi(x))$,
$v_L(x) = g (\psi(x))$ at nodes of level $L$ converge as $L\to \infty$ to $(u,v)$ the unique solution to
\eqref{ED1.general.9977} with \eqref{ED2.general.9977}
in the whole tree. 
\end{theorem}

\begin{proof}
From the estimates that we have proved in the previous section for the unique
solution $(u,v)$ to \eqref{ED1.general.9977} with \eqref{ED2.general.9977}
in the whole tree we know that, given $\eta>0$ there exists $L$ large enough
such that
we have
$$
u(x) \leq \max_{I_x} f + \eta, \qquad \mbox{and} \qquad v(x) \leq \max_{I_x} g + \eta
$$
for every $x$ at level $L$. 

On the other hand, since $f$ and $g$ are continuous, it holds that 
$$
|u_L(x) - \max_{I_x} f | = |f(\psi(x)) - \max_{I_x} f | < \eta, \quad \mbox{and} \quad 
|v_L(x) - \max_{I_x} g | = |g(\psi(x)) - \max_{I_x} g |< \eta
$$
for every $x$ at level $L$ with $L$ large enough. 

Therefore,
$(u,v)$ and $(u_L,v_L)$ are two solutions to the system \eqref{ED1.general.9977} 
in the finite subgraph of the tree with nodes of level less than $L$ that verify
$$
u(x) < u_L(x)  + 2 \eta, \qquad \mbox{and} \qquad v(x) < v_L(x)  + 2\eta
$$
for every $x$ at level $L$.
Now, since $(u_L(x)  + 2 \eta, v_L(x)  + 2 \eta)$ and $(u,v)$ are two solutions 
to \eqref{ED1.general.9977} in the subgraph of the tree with nodes of level less than $L$
that are ordered at its boundary (the set of nodes of level $L$) and the comparison 
principle can be used in this context we conclude that 
$$
u(x) \leq \liminf_{L\to \infty} u_L(x) , \qquad \mbox{and} \qquad v(x) \leq \liminf_{L\to \infty}  v_L(x).
$$

A similar argument starting with
$$
u(x) \geq \min_{I_x} f + \eta, \qquad \mbox{and} \qquad v(x) \geq \min_{I_x} g + \eta
$$ 
for every $x$ at level $L$ with $L$ large gives
$$
u(x) \geq \limsup_{L\to \infty} u_L(x) , \qquad \mbox{and} \qquad v(x) \leq \limsup_{L\to \infty}  v_L(x)
$$
and completes the proof. 
\end{proof}

{\bf Acknowledgements.} 

C. Mosquera is partially supported by grants UBACyT 20020170100430BA (Argentina), PICT 2018--03399 (Argentina) and PICT 2018--04027 (Argentina). A. Miranda and J. Rossi are partially supported by grants
CONICET grant PIP GI No 11220150100036CO
(Argentina), PICT-2018-03183 (Argentina) and UBACyT grant 20020160100155BA (Argentina).


\end{document}